\documentclass{amsproc}

\newtheorem{theorem}{Theorem}[section]
\newtheorem{lemma}[theorem]{Lemma}
\newtheorem{corollary}[theorem]{Corollary}
\newtheorem{proposition}[theorem]{Proposition}
\newtheorem{claim}{Claim}

\theoremstyle{definition}
\newtheorem{definition}[theorem]{Definition}
\newtheorem{example}[theorem]{Example}

\theoremstyle{remark}
\newtheorem{remark}[theorem]{Remark}

\numberwithin{equation}{section}

\begin{document}

\title[Remez Type Inequalities and Morrey-Campanato Spaces]{Remez Type Inequalities and Morrey-Campanato Spaces on Ahlfors Regular Sets}

\author{Alexander Brudnyi}
\address{Department of Mathematics and Statistics, University of Calgary,
Calgary, Canada}
\email{albru@math.ucalgary.ca}
\thanks{The first author was supported in part by NSERC}

\author{Yuri Brudnyi}
\address{Department of Mathematics, Technion,
Haifa, Israel} \email{ybrudnyi@math.technion.ac.il}


\dedicatory{This paper is dedicated to our friend Michael Cwikel
with respect and sympathy.}


\begin{abstract}
The paper presents several new results on Remez type inequalities
for real and complex polynomials in $n$ variables on Ahlfors regular
subsets of Lebesgue $n$-measure zero. As an application we prove an
extension theorem for Morrey-Campanato spaces defined on such sets.
\end{abstract}

\maketitle
\section{Introduction}
\font\mbn=msbm10 scaled \magstep1 \font\mbs=msbm7 scaled \magstep1
\font\mbss=msbm5 scaled \magstep1
\newfam\mbff
\textfont\mbff=\mbn \scriptfont\mbff=\mbs
\scriptscriptfont\mbff=\mbss\def\mbf{\fam\mbff}
\def\Re{{\mbf R}}
\def\Q{{\mbf Q}}
\def\Z{{\mbf Z}}
\def\Co{{\mbf C}}
\def\To{{\mbf T}}
\def\P{{\mbf P}}
\def\Di{{\mbf D}}
\def\Bo{{\mbf B}}
\def\F{{\mbf F}}
\def\N{{\mbf N}}
\def\H{{\mbf H}}
Recently there has been a considerable interest in multi-dimensional
analogs of the classical Remez polynomial inequality in connection
with various problems of Analysis, see, e.g., surveys [BB] and [G]
and references therein. The first result of this type is the
Yu.Brudnyi-Ganzburg inequality [BG]. In its formulation, we let
$V\subset\mathbf{R}^{n}$ be a convex body, $\omega\subset V$ a
measurable subset, and $\mathcal{L}_{n}$ be the Lebesgue measure on
$\mathbf{R}^{n}$.

Suppose that $p\in\mathbf{R}[x_{1},\dots, x_{n}]$ is a real
polynomial of degree $k$ on $\mathbf{R}^{n}$. Then the following
inequality holds:
\begin{equation}\label{e1}
\sup_{V}|p|\leq
T_{k}\left(\frac{1+\beta_{n}(\lambda)}{1-\beta_{n}(\lambda)}\right)\sup_{\omega}|p|.
\end{equation}
Here $T_{k}$ is the Chebyshev polynomial of degree $k$,
$\beta_{n}(\lambda):=(1-\lambda)^{1/n}$ and $\lambda:=\frac{\mathcal
{L}_{n}(\omega)}{\mathcal{L}_{n}(V)}$.

This inequality is sharp and for $n=1$ coincides with the classical
Remez inequality.

Most of the known applications of this inequality in Analysis use
the following corollary of (\ref{e1}):
\begin{equation}\label{e2}
\sup_{V}|p|\leq\left(\frac{4n}{\lambda}\right)^{k}\sup_{\omega}|p|
\end{equation}
which is easily derived from (\ref{e1}) using properties of the
Chebyshev polynomials.

In an actively developing field of modern mathematics, analysis on
fractal sets, see, e.g., [Tr] and references therein, one requires a
generalization of inequality (\ref{e2}) for fractal sets. In such a
generalization $\omega$ is a subset of Lebesgue measure $0$ in a
Euclidean ball $B\subset\mathbf{R}^{n}$. Since zero sets of real
polynomials on $\mathbf{R}^{n}$ have Hausdorff dimension $\leq n-1$,
to obtain a finite bound for $\sup_{B}|p|/\sup_{\omega}|p|$ one
assumes also that the Hausdorff dimension of $\omega$ is more than $n-1$.
Further, it is natural to estimate the above ratio by a function
depending on the Hausdorff measures of $B$ and $\omega$. Specifically,
let $\mathcal{H}_{s}$ denote the $s$-Hausdorff measure on
$\mathbf{R}^{n}$, $0<s\leq n$; in particular, $\mathcal{H}_{n}$
coincides with $\mathcal{L}_{n}$ up to a factor depending only on
$n$. In the present paper we study Remez type inequalities of the
following form
\begin{equation}\label{re3}
\sup_{B}|p|\leq\phi(\lambda)\sup_{\omega}|p|,
\end{equation}
where $p$ is a real polynomial on $\mathbf{R}^{n}$ or a holomorphic
polynomial on $\Co^{n}$, $B$ is a Euclidean ball in $\mathbf{R}^{n}$
or $\Co^{n}$, respectively, and $\omega\subset B$ is a subset of 
finite Hausdorff $s$-measure with $n-1<s\leq n$ in the real case and
$2n-2<s\leq 2n$ in the complex one. Also,
$$
\lambda:=\frac{\{\mathcal{H}_{s}(\omega)\}^{m/s}}{\mathcal{
H}_{m}(B)},
$$
where $m=n$ in the real case and $m=2n$ in the complex case.

 For many applications (related, e.g., to reverse
H\"{o}lder inequalities or BMO-properties of functions) it is
crucial that $\phi$ in (\ref{re3}) is a power function in $\lambda$.
Inequalities of the form (\ref{re3}) with such a function will be
referred to as {\em strong Remez type inequalities}. However, in
applications related to trace and extension theorems for classical
spaces of differentiable functions, see, in particular, [Ca], [YBr],
[J], it suffices to use inequalities of the form (\ref{re3}) with a
function $\phi$ whose dependence of $\lambda$ is not specified. In
this case the only required information is the monotonicity of $\phi$ in
$\lambda$. Such inequalities will be referred to as {\em weak Remez
type inequalities}.

The existence of inequalities (\ref{re3}) for $n=1$ was first
demonstrated in [ABr] where strong Remez type inequalities were
proved for ({\em Ahlfors}) $s$-regular sets $\omega$ in $\mathbf{R}$
or $\Co$ with $0<s\leq 1$ for real $p$ and with $0<s\leq 2$ for
holomorphic ones. Moreover, it was proved in [ABr, Prop.$\ \!$ 3]
that $s$-regularity is necessary for the validity of such an inequality.

Let us recall the definition of Ahlfors regular sets.

For a subset $K\subset\mathbf{R}^{n}$ and a point $x\in K$ by
$B_{r}(x;K)$ we denote the intersection with $K$ of an open
Euclidean ball in $\mathbf{R}^{n}$ centered at $x$ of radius $r$.
\begin{definition}\label{d1}
A subset $K\subset\mathbf{R}^{n}$ is said to be (Ahlfors)
$s$-regular if there is a positive number $a$ such that for every
$x\in K$ and $0<r\leq diam(K)$
\begin{equation}\label{e3}
\mathcal{H}_{s}(B_{r}(x;K))\leq ar^{s}.
\end{equation}
\end{definition}

The class of these sets will be denoted by $\mathcal{A}_{n}(s,a)$.
\begin{definition}\label{d2}
A subset $K\subset \mathcal{A}_{n}(s,a)$ is said to be an $s$-set if
there is a positive number $b$ such that for every $x\in K$ and
$0<r\leq diam(K)$
\begin{equation}\label{e4}
br^{s}\leq \mathcal{H}_{s}(B_{r}(x;K)).
\end{equation}
\end{definition}
We denote this class by $\mathcal{A}_{n}(s,a,b)$.

The class of $s$-sets, in particular, contains compact Lipschitz
$s$-manifolds (with integer $s$), Cantor type sets and self-similar
sets (with arbitrary $s$), see, e.g., [JW, p.\ $\!$29] and [Ma,
Sect.\ $\!$4.13].

In the present paper we establish inequalities of form (\ref{re3})
for $s$-regular sets $\omega\in \mathcal{A}_{n}(s,a)$ with $\phi$
depending also on $s$, $n$, $k:=deg\ \!p$ and $a$. We prove strong
Remez type inequalities for holomorphic polynomials using a
technique of Algebraic Geometry. For the real case, strong Remez type
inequalities are true for dimensions $n=1,2$ but the problem is open
for $n>2$. On the other hand, weak Remez type inequalities are
valid in this case, the proof is outlined in [BB]. For the convenience
of the reader we present this proof below.

In the final section of the paper we present an extension theorem for
Morrey-Campanato spaces defined on $s$-sets with $n-1<s\leq n$. For
$s=n$ this result is proved in [YBr].

Observe that a weak (or strong) Remez type inequality implies the,
so-called, Markov inequality for $s$-regular sets stating that for
some $c=c(F,n,k)$
\begin{equation}\label{markov}
\max_{F\cap B}|\nabla p|\leq\frac{c}{r}\max_{F\cap B}|p|
\end{equation}
where $F\in \mathcal{A}_{n}(s,a)$, $n-1<s\leq n$, $B$ is a closed
Euclidean ball of radius $r$ centered at $F$ and
$p\in\mathbf{R}[x_{1},\dots, x_{n}]$ is a polynomial of degree $k$.
For $s$-sets this result was established by other methods in [JW,
Sec.\ $\!$II.1.3].
\section{Formulation of Main Results}
\subsection{} We start with strong Remez type inequalities for
holomorphic polynomials on $\Co^{n}$.

Let $X\subset\Co^{n}$ belong to $\mathcal{A}_{2n}(s,a)$,
$s=2n-2+\alpha$, $\alpha>0$. Let $p$ be a holomorphic polynomial on
$\Co^{n}$ of degree $k$.
\begin{theorem}\label{remez1}
For any Euclidean ball $B\subset\Co^{n}$ and an
$\mathcal{H}_{s}$-measurable subset $\omega\subset X\cap B$ one has
$$
\sup_{B}|p|\leq\left(\frac{c_{1}\mathcal{H}_{2n}(B)}{\{\mathcal
{H}_{s}(\omega)\}^{2n/s}}\right)^{c_{2}k}\sup_{\omega}|p|
$$
where $c_{1}$ depends on $a$, $n$, $k$, $\alpha$ and $c_{2}>0$
depends on $\alpha$.
\end{theorem}
\begin{corollary}\label{remez2}
Let $X\in \mathcal{A}_{2n}(s,a,b)$. Let $B=B_{r}(x;X)$, $x\in X$,
$r>0$, and $\omega\subset B$ be $\mathcal{H}_{s}$-measurable. Then
for a holomorphic polynomial $p$ of degree $k$ the following is
true:
$$
\sup_{B}|p|\leq\left(\frac{c_{1}\mathcal{H}_{s}(B)}{\mathcal
{H}_{s}(\omega)}\right)^{c_{2}k}\sup_{\omega}|p|
$$
where $c_{1}$ depends on $a$, $b$, $n$, $k$, $\alpha$ and $c_{2}$
depends on $\alpha$.
\end{corollary}
\begin{corollary}\label{bmo}
Let $X\subset\Co^{n}$ be an $s$-set with $s$ as above. Then for any
holomorphic polynomial $p$ the function $\ln |p|\in
BMO(X,\mathcal{H}_{s})$.
\end{corollary}
Another corollary is the following reverse H\"{o}lder inequality.
\begin{corollary}\label{holder}
Under assumptions of Theorem \ref{remez1} for $1\leq l\leq\infty$
one has
$$
\left(\frac{1}{\mathcal{H}_{s}(B_{r}(x;X))}\int_{B_{r}(x;X)}|p|^{l}\
\!d\mathcal{H}_{s}\right)^{1/l}\leq
C\left(\frac{1}{\mathcal
{H}_{s}(B_{r}(x;X))}\int_{B_{r}(x;X)}|p|\ \!d\mathcal{H}_{s}\right)
$$
where $C$ depends on $k$, $n$, $\alpha$, $a$ and $b$.
\end{corollary}
\begin{remark}
{\rm The results of this subsection for $n=1$ were proved in [ABr]
with $c_{1}$ independent of $k$. An interesting open question is
whether $c_{1}$ is independent of $k$ in the general case, as well.}
\end{remark}
\subsection{} In this part we present a general form of
weak Remez type inequalities for real polynomials on $\Re^{n}$.
\begin{theorem}\label{weak}
Assume that $U\subset\Re^{n}$ is a bounded open set and
$\omega\subset U$ belongs to $\mathcal{A}_{n}(s,a)$ with $n-1<s\leq
n$. Assume also that
$$
\lambda:=\frac{\{\mathcal{H}_{s}(\omega)\}^{n/s}}{\mathcal{
H}_{n}(U)}>0.
$$
Then there is a constant $C>1$ such that for every polynomial
$p\in\Re[x_{1},\dots, x_{n}]$ of degree $k$
\begin{equation}\label{weak1}
\left(\frac{1}{\mathcal{H}_{n}(U)}\int_{U}|p|^{r}\
\!d\mathcal{H}_{n}\right)^{1/r}\leq C\left(\frac{1}{\mathcal
{H}_{s}(\omega)}\int_{\omega}|p|^{q}\
\!d\mathcal{H}_{s}\right)^{1/q}.
\end{equation}
Here $0<q,r\leq\infty$ and $C$ depends on $U$, $n$, $q$, $r$, $s$,
$k$, $a$ and $\lambda$ and is increasing in $1/\lambda$. In
particular, for $q=r=\infty$ we obtain the weak Remez type
inequality of the form (\ref{re3}).
\end{theorem}
\subsection{} Let $X\subset\Re^{n}$ be a measurable set of positive
Hausdorff $s$-measure. By $\mathcal{K}_{X}$ we denote the family of closed
cubes in $\Re^{n}$ with centers at $X$ and ``radii`` \penalty-10000
($:=\frac{1}{2}$ lengthside) at most $4 diam\ \!X$. We write $Q_{r}(x)$ for
the cube of radius $r$ and center $x$ and denote by $X_{r}(x)$ the
set $Q_{r}(x)\cap X$ for $x\in X$.

In order to introduce the basic concept, Morrey-Campanato space on
$X$, we denote by $L_{q}(X)$, $1\leq q\leq\infty$, the linear space
of $\mathcal{H}_{s}$-measurable functions on $X$ equipped with
norm
\begin{equation}\label{morrey1}
||f||_{q}:=\left(\int_{X}|f|^{q}\ \!d\mathcal{H}_{s}\right)^{1/q},\
\ \ 0\leq q\leq\infty,
\end{equation}
and use the following
\begin{definition}\label{morrey2}
The local best approximation of order $k\in\Z_{+}$ is a function
$\mathcal{E}_{k}:L_{q}(X)\times\mathcal{K}_{X}\to\Re_{+}$ given for
$Q=Q_{r}(x)$ by
\begin{equation}\label{morrey3}
\mathcal{E}_{k}(f;Q):=\inf_{p}\left\{\frac{1}{\mathcal{H}_{s}(X_{r}(x))}
\int_{X_{r}(x)}|f-p|^{q}\ \!d\mathcal{H}_{s}\right\}^{1/q}
\end{equation}
where $p$ runs over the space
$\mathcal{P}_{k-1}\subset\Re[x_{1},\dots,x_{n}]$ of polynomials of
degree $k-1$.
\end{definition}

For $k=0$ we let $\mathcal{P}_{k-1}:=\{0\}$; hence
$\mathcal{E}_{0}(f;Q)$ is the normalized $L_{q}$-norm of $f$ on
$X_{r}(x)$.

Let now $\omega:\Re_{+}\to\Re_{+}$ be a monotone function on
$\Re_{+}:=(0,\infty)$ (it may be a constant).
\begin{definition}\label{morrey4}
The (generalized) Morrey-Campanato space $\dot{C}_{q}^{k,\omega}(X)$
is defined by seminorm
$$
|f|_{\dot{C}_{q}^{k,\omega}(X)}:=\sup\left\{\frac{\mathcal{E}_{k}(f;Q)}
{\omega(r_{Q})}\ :\ Q\in\mathcal{K}_{X}\right\}
$$
where $r_{Q}$ denotes the radius of $Q$.
\end{definition}

For $X$ being a domain in $\Re^{n}$ and $s=n$ this space coincides with
the Morrey space $\mathcal{M}_{q}^{\lambda}$ [Mo] (for $k=0$,
$\omega(t)=t^{\lambda}$, $-n<\lambda<0$), the BMO-space [JN] (for $k=1$,
$\omega(t)=const$) and the Campanato space [Ca] (for $k\geq 1$,
$\omega(t)=t^{\lambda}$, $\lambda>0$).

To formulate the main result we also need
\begin{definition}\label{morrey5}
Let $\omega:\Re_{+}\to\Re_{+}$ be nondecreasing such that
$$
\omega(+0)=0\ \ \ {\rm and}\ \ \ t\to\frac{\omega(t)}{t^{k}}\ \ \
{\rm be\ nonincreasing}.
$$

The $\omega$ is said to be a quasipower $k$-majorant if
$$
C_{\omega}:=\sup_{t>0}\left\{
\frac{1}{\omega(t)}\int_{0}^{t}\frac{\omega(u)}{u}
\ \!du\right\}<\infty.
$$

The Lipschitz space $\dot{\Lambda}^{k,\omega}(\Re^{n})$ of order
$k\geq 1$ consists of locally bounded functions $f:\Re^{n}\to\Re$
such that the seminorm
$$
|f|_{\dot{\Lambda}^{k,\omega}(\Re^{n})}:=\sup\left\{\frac{|\Delta_{h}^{k}f(x)|}
{\omega(|h|)}\ :\ x,h\in\Re^{n}\right\}
$$
is finite.
\end{definition}

Here $|h|$ is the Euclidean norm of $h$ and
$$
\Delta_{h}^{k}f(x):=\sum_{j=0}^{k}(-1)^{k-j}{k \choose j}f(x+jh).
$$
\begin{example}\label{ex1}
{\rm  Choosing in this definition $\omega(t):=t^{\lambda}$, $0<\lambda\leq k$,
we obtain the (homogeneous) Besov space $B_{\infty}^{\lambda}(\Re^{n})$.
Let us recall that it coincides with the Sobolev space 
$\dot{W}_{\infty}^{k}(\Re^{n})$ for $\lambda=k$, the H\"{o}lder space 
$C^{l,\alpha}(\Re^{n})$ for $\lambda=l+\alpha$, $l$ is an integer and
$0<\alpha<1$, and with the Marchaud-Zygmund space for $\lambda$ integer and
$0<\lambda<k$. In the last case, the corresponding seminorm is
$$
|f|_{B_{\infty}^{\lambda}(\Re^{n})}:=\max_{|\alpha|=\lambda-1}\sup_{h}
\frac{||\Delta_{h}^{2}(D^{\alpha}f)||_{C(\Re^{n})}}{|h|}. 
$$
}
\end{example}
\begin{theorem}\label{ybr}
Let $X\subset\Re^{n}$ be an $s$-set with $n-1<s\leq n$ and
$\omega$ be a quasipower $k$-majorant.
Then there is a linear continuous extension operator
$T_{k}:\dot{C}_{q}^{k,\omega}(X)\to\dot{\Lambda}^{k,\omega}(\Re^{n})$.

In particular, $\dot{C}_{q}^{k,\omega}(X)$ is isomorphic
to the trace space
$\dot{\Lambda}^{k,\omega}(\Re^{n})|_{X}$.
\end{theorem}

For $\omega(t)=t^{\lambda}$, $0<\lambda\leq k$, and $n$-sets 
this result was proved in [YBr] in a different way.
\begin{remark}\label{ybr1}
{\rm (a) Theorem \ref{ybr} is also true for the nonhomogeneous
Morrey-Campanato space $C_{q}^{k,\omega}(X)$ defined by 
(quasi)-norm
$$
||f||_{C_{q}^{k,\omega}(X)}:=||f||_{q}+|f|_{\dot{C}_{q}^{k,\omega}(X)}
$$
The target space of the extension operator is now the Banach space
$\Lambda^{k,\omega}(\Re^{n})$ defined by norm
$$
||f||_{\Lambda^{k,\omega}(\Re^{n})}:=\sup_{Q_{0}}|f|+|f|_{\dot{\Lambda}^{k,\omega}(\Re^{n})}
$$
where $Q_{0}:=[0,1]^{n}$.\\
(b) For $0<q<1$ the extension operator exists also but it is only
nonlinear (homogeneous and bounded). \\
(c) Let $X_{j}\subset\Re^{n_{j}}$ be $s_{j}$-sets, 
$n_{j}-1< s_{j}\leq n_{j}$, $1\leq j\leq k$. Then under the same assumptions for $\omega$, Theorem \ref{ybr} is also valid for $X:=X_{1}\times\cdots X_{k}\subset\Re^{n}$, $n:=n_{1}+\dots +n_{k}$. Here the space $\dot{C}_{q}^{k,\omega}(X)$ is defined with respect to the tensor product of $\mathcal{H}_{s_{j}}$-measures. This can be proved similarly to the proof of Theorem \ref{ybr} based on the corresponding Remez type inequalities for such $X$. 
}
\end{remark}
\section{Holomorphic Polynomials}
\subsection{Complex Algebraic Varieties and $s$-sets.}
\subsubsection{}
In this section we gather some standard facts of Complex Algebraic
Geometry. For the background and the proofs see, e.g., books [M] and
[GH].

By $\Co\P^{n}$ we denote the $n$-dimensional complex projective
space with homogeneous coordinates $(z_{0}:\cdots : z_{n})$. The
complex vector space $\Co^{n}$ is a dense open subset of $\Co\P^{n}$
defined by $z_{0}\neq 0$. The hyperplane at $\infty$,
$H:=\{(z_{0}:\cdots :z_{n})\in\Co\P^{n}\ :\ z_{0}=0\}$, can be
naturally identified with $\Co\P^{n-1}$ and $\Co\P^{n}=\Co^{n}\cup
H$.

A closed subset $X\subset\Co^{n}$ defined as the set of zeros of a
family of holomorphic polynomials on $\Co^{n}$ is called an {\em
affine algebraic variety}. By $dim_{\Co}X$ we denote the (complex)
{\em dimension} of $X$, i.e., the maximum of complex dimensions of
complex tangent spaces at smooth points of $X$.

Assume that an affine algebraic variety $X\subset\Co^{n}$ has pure
dimension $k\geq 1$, i.e., dimensions of complex tangent spaces at
smooth points of $X$ are the same. Then its closure $\overline X$ in
$\Co\P^{n}$ is a projective variety of pure dimension $k$, and
$dim_{\Co}(H\cap\overline X)=k-1$.

Any linear subspace of dimension $n-k$ in $\Co\P^{n}$ meets
$\overline X$, but there is a linear subspace $L\subset H$ of
dimension $n-k-1$ such that $L\cap\overline X=\emptyset$. Moreover,
for a generic $(n-k)$-dimensional subspace of $\Co\P^{n}$ its
intersection with $\overline X$ consists of a finite number of
points. The number of these points is called the {\em degree} of
$\overline X$ and is denoted $deg\ \!\overline X$. For instance, if
$X$ as above is defined as the set of zeros of holomorphic
polynomials $p_{1},\dots, p_{n-k}$ on $\Co^{n}$ of degrees
$d_{1},\dots, d_{n-k}$, respectively, then by the famous Bezout
theorem $deg\ \!\overline{X}\leq d_{1}\cdots d_{n-k}$.

Let $L\subset H$ be a linear subspace of dimension $n-k-1$ which
does not intersect $\overline X$. This subspace defines a projection
$\phi_{L}:\Co\P^{n}\to\Co\P^{k}$ as follows.

Fix a linear subspace of dimension $k$ in $\Co\P^{n}$ disjoint from
$L$. We will simply call it $\Co\P^{k}$. If $w\in\Co\P^{n}\setminus
L$, then $w$ and $L$ span an $(n-k)$-dimensional linear subspace
which meets $\Co\P^{k}$ in a unique point $\phi_{L}(w)$. The map
$\phi_{L}$ sends $w$ to $\phi_{L}(w)$. Further,
$\Co^{n}\subset\Co\P^{n}\setminus L$, and, with a suitable choice of
linear coordinates, $\phi_{L}|_{\Co^{n}}:\Co^{n}\to\Co^{k}$ is the
standard projection: $(z_{1},\dots,z_{n})\mapsto (z_{1},\dots,
z_{k})$.

The map $\phi_{L}|_{X}:X\to\Co^{k}$ is a surjection and is a
branched covering over $\Co^{k}$ whose order $\mu$, i.e., the number
of points $\phi_{L}^{-1}(y)\cap X$ for a generic $y\in\Co^{k}$, is
$deg\ \!\overline X$. Then $X$ is a complex subvariety of a pure
$k$-dimensional algebraic variety $\widetilde X$ defined as the set
of zeros of holomorphic polynomials $p_{i}$, $1\leq i\leq n-k$, of
the form
\begin{equation}\label{form}
p_{i}(z_{1},\dots, z_{n})=z_{k+i}^{\mu}+\sum_{1\leq l\leq
\mu}b_{il}(z_{1},\dots, z_{k})z_{k+i}^{\mu-l}
\end{equation}
where $b_{il}$ is a holomorphic polynomial of degree $\leq l$ on
$\Co^{k}$. Moreover, let $S\subset\Co^{k}$ be the branch locus of
$\phi_{L}|_{X}$. If $w\in\Co^{k}\setminus S$, then $b_{il}(w)$ is
the $l$-th elementary symmetric function in $z_{k+i}(w^{(1)}),\dots,
z_{k+i}(w^{(\mu)})$, where $\phi_{L}^{-1}(w)\cap X=(w^{(1)},\dots,
w^{(\mu)})$. (Recall that the elementary symmetric functions $s_{i}$
in $\xi_{1},\dots,\xi_{\mu}$ are defined from the identity
$\prod_{1\leq l\leq
\mu}(t-\xi_{l})=t^{\mu}+s_{1}t^{\mu-1}+\cdots+s_{\mu}$ of
polynomials in variable $t$.) Since $dim_{\Co}X=dim_{\Co}\widetilde
X=k$, $X$ is the union of some irreducible components of $\widetilde
X$.

Next, the Fubini-Studi metric on $\Co\P^{n}$ is a Riemannian
metric defined by the associated $(1,1)$-form
$\omega:=\frac{\sqrt{-1}}{2\pi}\partial\overline\partial\ln(|z_{0}|^{2}+\cdots
|z_{n}|^{2})$, $(z_{0}:\cdots : z_{n})\in\Co\P^{n}$. For a
$k$-dimensional projective variety $\overline X$ as above the
$(k,k)$-form $\wedge^{k}\omega$ determines a Borel measure
$\mu_{\overline{X}}$ on $\overline X$, 
\begin{equation}\label{meas1}
\mu_{\overline{X}}(U):=\int_{U}\wedge^{k}\omega
\end{equation}
where $U\subset\overline X$ is a Borel subset. Moreover,
\begin{equation}\label{meas1'}
\mu_{\overline{X}}(\overline{X})=deg\ \!\overline{X},
\end{equation}
see, e.g., [GH, Ch.\ $\!$1.5].

Let $\omega_{e}:=\frac{\sqrt{-1}}{2}\sum_{1\leq i\leq n}dz_{i}\wedge
d\overline{z}_{i}$ be the Euclidean K\"{a}hler form determining the
Euclidean metric on $\Co^{n}$. Then $\omega$ and $\omega_{e}$ are
equivalent on every compact subset $K\subset\Co^{n}$ where the
constants of equivalence depend on $K$ and $n$ only. In particular,
the Fubini-Studi  and the Euclidean metrics, and the $(k,k)$-forms
$\wedge^{k}\omega_{e}$ and $\wedge^{k}\omega$ are equivalent on
every such $K$. Let $\mu_{e,X}$ be a Borel measure on a pure
$k$-dimensional affine algebraic variety $X$ defined by the formula
\begin{equation}\label{meas2}
\mu_{e,X}(U):=\int_{U}\wedge^{k}\omega_{e}
\end{equation}
where $U\subset X$ is a Borel subset. Then for every compact subset
$K\subset\Co^{n}$ the measures $\mu_{\overline{X}}|_{K\cap X}$ and
$\mu_{e,X}|_{K\cap X}$ are equivalent with the constants of
equivalence depending on $K$, $k$ and $n$ only.

\subsubsection{} In this section we establish a relation between
complex algebraic varieties and $s$-sets.
\begin{theorem}\label{dset}
Let $X\subset\Co^{n}$ be an affine algebraic variety of pure
dimension $k\geq 1$ such that $deg\ \!\overline X\leq\mu$. Then
$X\in \mathcal{A}_{2n}(2k,a,b)$ where $a$ and $b$ depend on $k$,
$\mu$ and $n$ only.
\end{theorem}
\begin{proof} We will prove that
\begin{equation}\label{dset0}
br^{2k}\leq \mu_{e,X}(B_{r}(x;X))\leq ar^{2k}
\end{equation}
with $a$ and $b$ depending on $k$, $\mu$ and $n$ only, where $X$
satisfies the assumptions of the theorem, $x\in X$ and $\mu_{e,X}$
is the measure on $X$ determined in (\ref{meas2}). From here
applying [JW, Sec. II.1.2, Th. 1] we get the desired statement.

Since $deg\ \!\overline{X}=deg\ \!\overline{x+X}$ and
$\mu_{e,X}(U)=\mu_{e,x+X}(x+U)$ for all $x\in\Co^{n}$ and all Borel
subsets $U\subset X$, without loss of generality we may assume that
$0\in X$ and prove (\ref{dset0}) for $B_{r}(0;X)$ only. Since
$\mu_{e,\lambda X}(\lambda U)=\lambda^{2k}\mu_{e,X}(U)$ and $deg\
\!\overline{\lambda X}=deg\ \!\overline{X}$ for $\lambda>0$,
$x\in\Co^{n}$, and Borel subsets $U\subset\Co^{n}$, it suffices to
prove that
\begin{equation}\label{dset1}
b\leq \mu_{e,X}(B_{1}(0;X))\leq a
\end{equation}
where $a$ and $b$ depend on $k$, $\mu$ and $n$ only.

First we will prove the left-side inequality in (\ref{dset1}). Let
$\{X_{l}\}_{l\in\N}$ be a sequence of affine algebraic varieties
containing $0$ and satisfying the hypotheses of the theorem such
that

\begin{equation}\label{dset2}
\inf_{X}\mu_{e,X}(B_{1}(0;X))=\lim_{l\to\infty}
\mu_{e,X_{l}}(B_{1}(0;X_{l})).
\end{equation}

Here the infimum is taken over all $X$ containing $0$ and satisfying
the conditions of the theorem. Consider the sequence
$\{\overline{X}_{l}\}_{l\in\N}$ of pure $k$-dimensional projective
subvarieties of $\Co\P^{n}$. Since $\Co\P^{n}$ is a compact
manifold, one can choose a subsequence of
$\{\overline{X}_{l}\}_{l\in\N}$ converging in the Hausdorff metric
defined on compact subsets of $\Co\P^{n}$ to a compact set, say,
$Y$. Without loss of generality we may assume that
$\{\overline{X}_{l}\}_{l\in\N}$ itself converges to $Y$.
\begin{lemma}\label{dset3}
There are a linear subspace $L\subset\Co\P^{n}$ of dimension $n-k-1$
and a number $N\in\N$ such that $L\cap(\{\overline{X}_{l}\}_{l\geq
N}\cup Y)=\emptyset$.
\end{lemma}
\begin{proof} We prove the result by induction on $n-k$, the codimension of
$\overline{X}_{l}$ in $\Co\P^{n}$.

For $n-k=1$ every $\overline{X}_{l}$ being a projective hypersurface
of degree $\leq\mu$ is defined as the set of zeros of a holomorphic
homogeneous polynomial $p_{l}$ of degree $\leq\mu$:
$$
\overline{X}_{l}:=\{(z_{0}:\cdots : z_{n})\in\Co\P^{n}\ :\
p_{l}(z_{0},\dots,z_{n})=0\}.
$$
Without loss of generality we may assume that $l_{2}$-norms of
vectors of coefficients of all $p_{l}$ are $1$. Then we can
choose a subsequence $\{p_{l_{s}}\}_{s\in\N}$ that converges
uniformly on compact subsets of $\Co^{n+1}$ to a nontrivial
(holomorphic) homogeneous polynomial $p$ of $deg\ \!p\leq\mu$. Next,
if $y\in Y$, then by the definition of the Hausdorff convergence
there is a sequence of points $\{x_{l}\}_{l\in\N}$, $x_{l}\in
X_{l}$, such that $\lim_{l\to\infty}x_{l}=y$. In particular, if
$y=(y_{0}:\cdots :y_{n})$ and $x_{l}=(x_{0l}:\cdots :x_{nl})$ with
$\max_{0\leq i\leq n}|y_{i}|\leq 1$, $\max_{0\leq i\leq
n}|x_{il}|\leq 1$, $l\in\N$, then
$$
p(y_{0},\dots,y_{n})=\lim_{s\to\infty}p_{l_{s}}(x_{0l_{s}},\dots,x_{nl_{s}})=0.
$$
Since $p\neq 0$, the latter implies that $Y$ belongs to a projective
hypersurface in $\Co\P^{n}$. In particular, $Y$ is nowhere dense in
$\Co\P^{n}$. Thus there is $z\in\Co\P^{n}\setminus Y$. And so there
is a neighbourhood $U$ of $Y$ in $\Co\P^{n}$ which does not contain
$z$. By the definition of the Hausdorff convergence this implies
that there is a number $N\in\N$ such that
$\{\overline{X}_{l}\}_{l\geq N}\subset U$ completing the proof of
the lemma for $n-k=1$.

Suppose now that the result is proved for $n-k>1$ and prove it for
$n-k+1$.

Since every $\overline{X}_{l}$ is contained in a projective
hypersurface in $\Co\P^{n}$ of degree $\leq\mu$ (see section 3.1.1),
by the induction hypothesis there are a number $N'\in\N$ and a point
$y\in\Co\P^{n}$ such that $y\not\in\{\overline{X}_{l}\}_{l\geq
N'}\cup Y$. The point $y$ determines a projection
$\phi_{y}:\Co\P^{n}\setminus\{y\}\to\Co\P^{n-1}$ as described in
section 3.1.1 (with $L:=\{y\}$). Set
$X_{l}'=\phi_{y}(\overline{X}_{l})$, $l\geq N'$, and
$Y'=\phi_{y}(Y)$. By the proper map theorem (see, e.g., [GH, Ch.\
$\!$0.2]) and the Chow theorem (see, e.g., [GH, Ch.\ $\!$1.3])
$X_{l}'$ are projective subvarieties of $\Co\P^{n-1}$. Also, by the
construction, cf. section 3.1.1, $dim_{\Co}X_{l}'=dim_{\Co}X_{l}$
and $deg\ \! X_{l}'\leq\mu$ for all $l\geq N'$. Moreover,
$\{X_{l}'\}_{l\geq N'}$ converges in the Hausdorff metric defined on
compact subsets of $\Co\P^{n-1}$ to $Y'$, because $\phi_{y}$ is
continuous in a neighbourhood of $\{\overline{X}_{l}\}_{l\geq
N'}\cup Y$. Since the codimension of $X_{l}'$ in $\Co\P^{n-1}$ is
$n-k$, by the induction hypothesis there are an integer number
$N\geq N'$ and a linear subspace $L'\subset\Co\P^{n-1}$ of dimension
$n-k-1$ which does not intersect $\{X_{l}'\}_{l\geq N}\cup Y'$. Then
$L=\phi_{y}^{-1}(L')\cup\{y\}$ is a linear subspace of $\Co\P^{n}$
of dimension $n-k$ which does not intersect
$\{\overline{X}_{l}\}_{l\geq N}\cup Y$.

This completes the proof of the lemma.
\end{proof}

Further, since $0\in\{X_{l}\}_{l\in\N}\cup Y$, there is a closed
Euclidean ball $\overline{B}_{r_{0}}(0)\subset\Co^{n}$ centered at
$0$ of radius $0<r_{0}\leq 1$ which does not intersect the $L$ of
the above lemma. Clearly,
$$
\mu_{e,X_{l}}(B_{1}(0;X_{l}))\geq
\mu_{e,X_{l}}(B_{r_{0}}(0;X_{l})),\ \ \ l\in\N.
$$
(As before, $B_{r_{0}}(0;X_{l}):=B_{r_{0}}(0)\cap X_{l}$.)
 Therefore to prove the left-side inequality in
(\ref{dset1}) it suffices to check that
\begin{equation}\label{au1}
\liminf_{l\to\infty}\mu_{e,X_{l}}(B_{r_{0}}(0;X_{l}))>0.
\end{equation}
Recall that the Fubini-Studi metric is equivalent to the Euclidean
metric on every compact subset $K\subset\Co^{n}$ with the constants
of equivalence depending on $K$ and $n$ only. Therefore there is a
closed ball $B$ in the Fubini-Studi metric centered at $0$ and of
radius $s_{0}>0$ depending on $r_{0}$ and $n$ only such that
$B\subset\overline{B}_{r_{0}}(0)$.  Since $\mu_{e,X_{l}}$ is
equivalent to $\mu_{\overline{X}_{l}}$ on
$\overline{B}_{r_{0}}(0;X_{l})$ with the constants of equivalence
depending on $r_{0}$, $k$ and $n$ only, see section 3.1.1,
inequality (\ref{au1}) follows from the inequality
\begin{equation}\label{au2}
\liminf_{l\to\infty}\mu_{\overline{X}_{l}}(B\cap X_{l})>0.
\end{equation}

Let us check the last inequality. Diminishing, if necessary, $r_{0}$
we can find a hyperplane $L'\subset\Co\P^{n}$ which contains $L$
from Lemma \ref{dset3} and does not intersect
$\overline{B}_{r_{0}}(0)$. Let $T:\Co^{n+1}\to\Co^{n+1}$ be a
unitary transformation which induces an isometry $\overline
T:\Co\P^{n}\to\Co\P^{n}$ sending $L'$ to the hyperplane at $\infty$,
$H$. Then $\overline T(B)$ is a closed ball (in the Fubini-Studi
metric) in $\Co^{n}=\Co\P^{n}\setminus H$. By the definition of $T$,
$deg\ \! T(\overline{X}_{l})=deg\ \!\overline{X}_{l}$ and
$\mu_{\overline{X}_{l}}(U)=\mu_{T(\overline{X}_{l})}(T(U))$ for a
Borel subset $U\in\overline{X}_{l}$. These facts and the above
equivalence of $\mu_{e,X_{l}}$ and $\mu_{\overline{X}_{l}}$ on
compact subsets of $\Co^{n}$ show that in the proof of (\ref{au2})
without loss of generality we may assume that $L'=H$.

Now, consider the projection $\phi_{L}:\Co^{n}\to\Co^{k}$ determined
as in section 3.1.1. Choosing suitable coordinates on $\Co^{n}$
we may and will assume that $\phi_{L}$ coincides with the projection
$(z_{1},\dots, z_{n})\mapsto (z_{1},\dots, z_{k})$. Then
$X_{l}:=\overline{X}_{l}\setminus H$ are algebraic subvarieties of
algebraic varieties $\widetilde X_{l}$ defined as sets of zeros of
families of polynomials $p_{il}$, $1\leq i\leq n-k$, $l\geq N$, of the
form (\ref{form}). Moreover, since
$L\cap(\{\overline{X}_{l}\}_{l\geq N}\cup Y)=\emptyset$, the
definition of $p_{il}$, see section 3.1.1, shows that for every $i$
polynomials $p_{il}$, $l\geq N$, are uniformly bounded on compact
subsets of $\Co^{n}$. Since $deg\ \!p_{il}\leq\mu$, we can find a
subsequence $\{l_{s}\}_{s\in\N}\subset\N$ such that
$\{p_{il_{s}}\}_{s\in\N}$ converge uniformly on compact subsets of
$\Co^{n}$ to polynomials $p_{i}$, $deg\ \!p_{i}\leq\mu$, of the form
(\ref{form}), $1\leq i\leq n-k$, and
$$
\lim_{s\to\infty}\mu_{\overline{X}_{l_{s}}}(B\cap
X_{l_{s}})=\liminf_{l\to\infty}\mu_{\overline{X}_{l}}(B\cap X_{l}).
$$
This implies easily that $Y\cap\Co^{n}$ with $Y$ from Lemma
\ref{dset3} is contained in the pure $k$-dimensional algebraic
variety $\widetilde Y$ defined as the set of zeros of polynomials
$p_{i}$, $1\leq i\leq n-k$.

In what follows by $\Delta_{r}^{l}:=\{(z_{1},\dots,
z_{l})\in\Co^{l}\ :\ \max_{1\leq i\leq l}|z_{i}|<r\}$ we denote the
open polydisk in $\Co^{l}$ centered at $0$ of radius $r$.

Since, by the definition, $\widetilde Y$ is a finite branched
covering over $\Co^{k}$ and $0\in \widetilde Y$, there is a polydisk
$\Delta_{\epsilon}^{n}=\Delta_{\epsilon}^{n-k}\times\Delta_{\epsilon}^{k}$
such that $\Delta_{\epsilon}^{n}\cap \widetilde Y\subset B\cap
\widetilde Y$ and $\phi_{L}:\Delta_{\epsilon}^{n}\cap \widetilde
Y\to\Delta_{\epsilon}^{k}$ is a finite branched covering over
$\Delta_{\epsilon}^{k}$ (for similar arguments see, e.g., the proof
of the preparatory Weierstrass theorem in [GH, Ch.\ $\!$0.1]). From
here using the fact that $\{p_{il_{s}}\}$ converges uniformly on
compact subsets of $\Co^{n}$ to $p_{i}$ for all $i$ and diminishing,
if necessary, $\epsilon$ we obtain analogously that there is a
number $N_{0}\in\N$ such that $\Delta_{\epsilon}^{n}\cap \widetilde
X_{l_{s}}\subset B\cap \widetilde X_{l_{s}}$ and
$\phi_{L}:\Delta_{\epsilon}^{n}\cap \widetilde
X_{l_{s}}\to\Delta_{\epsilon}^{k}$ are finite branched coverings
over $\Delta_{\epsilon}^{k}$ for all $s\geq N_{0}$. But
$\Delta_{\epsilon}^{n}\cap X_{l_{s}}$ is a (closed) complex
subvariety of $\Delta_{\epsilon}^{n}\cap \widetilde X_{l_{s}}$, and
$\phi_{L}(\Delta_{\epsilon}^{n}\cap X_{l_{s}})$ is an open subset of
$\Delta_{\epsilon}^{k}$ (because $0\in X_{l_{s}}$ and
$\phi_{L}:X_{l_{s}}\to\Co^{k}$ is a finite branched covering). Thus
by the proper map theorem, $\phi_{L}(\Delta_{\epsilon}^{n}\cap
X_{l_{s}})=\Delta_{\epsilon}^{k}$. (Here we used the fact that the
map $\phi_{L}:\Delta_{\epsilon}^{n}\cap
X_{l_{s}}\to\Delta_{\epsilon}^{k}$ is proper, because
$\phi_{L}:\Delta_{\epsilon}^{n}\cap \widetilde
X_{l_{s}}\to\Delta_{\epsilon}^{k}$ is proper and $X_{l_{s}}\cap
\Delta_{\epsilon}^{n}$ is a complex subvariety of $\widetilde
X_{l_{s}}\cap\Delta_{\epsilon}^{n}$.)

Let $\mathcal{L}_{2k}$ be the Lebesgue measure on $\Co^{k}$. Then by
the definition of $\mu_{\overline{X}_{l_{s}}}$ there is a constant
$c>0$ depending on $\mu$, $k$ and $n$ only such that
$$
\mu_{\overline{X}_{l_{s}}}(B\cap X_{l_{s}})\geq c \mathcal
{L}_{2k}(\phi_{L}(B\cap X_{l_{s}})).
$$
But for $s\geq N_{0}$ we have
$$
\mathcal{L}_{2k}(\phi_{L}(B\cap X_{l_{s}}))\geq \mathcal{
L}_{2k}(\Delta_{\epsilon}^{k})=\pi^{k}\epsilon^{2k}>0.
$$
The combination of the last two inequalities completes the proof of
(\ref{au2}) and thus the proof of  the left-side inequality in
(\ref{dset1}).

The right-side inequality in (\ref{dset1}) is obtained as follows,
see (\ref{meas1'}),
$$
\mu_{e,X}(B_{1}(0;X))\leq c(n,k)\mu_{\overline{X}}(B_{1}(0;X))\leq
c(n,k)\mu_{\overline{X}}(\overline{X})=c(n,k)deg\ \!\overline{X}\leq
c(n,k)\mu.
$$

The proof of Theorem \ref{dset} is complete.
\end{proof}
\subsection{Strong Remez Type Inequalities}
\subsubsection{Covering Lemmas} Our proof is based on a deep
generalization of the classical Cartan Lemma [C] discovered by Gorin
[GK]. We present a more general version of this result.

Let $X$ be a pseudometric space with pseudometric $d$. By
$\mathcal{F}:=\{\overline{B}_{r}(x)\subset X\ :\  d(x,y)\leq r,\
x,y\in X, r\geq 0\}$ we denote the set of closed balls in $X$. Let
$\xi:\mathcal{F}\to\mathbf{R}_{+}$ be a function satisfying the
following two properties:
\begin{enumerate}
\item
$$
\xi(\overline{B}_{r'}(x))\leq \xi(\overline{B}_{r''}(x))\ \ \ {\rm
for\ all}\ \ \ x\in X,\ r'\leq r''.
$$
\item
There is a numerical constant $A$ such that for any collection of
mutually disjoint balls $\{B_{i}\}\subset \mathcal{F}$,
$$
\sum_{i\geq 1}\xi(B_{i})\leq A
$$
\end{enumerate}

Consider a continuous strictly increasing nonnegative function
$\phi$ on $[0,\infty)$, $\phi(0)=0$, $\lim_{t\to\infty}\phi(t)>A$
which will be called a {\em majorant}.

For each point $x\in X$ we set $\tau(x)=\sup\{t\ :\ \xi(\overline
{B}_{t}(x))\geq\phi(t)\}$. It is easy to see that $\xi(\overline
{B}_{\tau(x)}(x))=\phi(\tau(x))$ and
$\sup_{x}\tau(x)\leq\phi^{-1}(A)<\infty$.

A point $x\in X$ is said to be {\em regular} (with respect to $\xi$
and $\phi$) if $\tau(x)=0$, i.e., $\xi(\overline{B}_{t}(x))<\phi(t)$
for all $t>0$. The next result shows that the set of regular points
is sufficiently large for an arbitrary majorant $\phi$.
\begin{lemma}\label{gorin}
Fix $\gamma\in (0,1/2)$. There is a sequence of balls
$B_{k}=\overline{B}_{t_{k}}(x_{k})$, $k=1,2,\dots$, which
collectively cover all irregular points such that $\sum_{k\geq
1}\phi(\gamma t_{k})<A$ (i.e., $t_{k}\to 0$).
\end{lemma}
The proof of this lemma for $\xi$ being a finite Borel measure on a
metric space $X$ is given by Gorin [GK]. His argument works also in
the general case.
\begin{proof} Let $0<\alpha<1$, $\beta>2$ be such that $\gamma<
\alpha/\beta$. We set $B_{0}=\emptyset$ and assume that the balls
$B_{0},\dots, B_{k-1}$ have been constructed. If
$\tau_{k}=\sup\{\tau(x)\ :\ x\not\in B_{0}\cup\dots\cup B_{k-1}\}$,
then there exists a point $x_{k}\not\in B_{0}\cup\dots B_{k-1}$ such
that $\tau(x_{k})\geq\alpha\tau_{k}$. We set $t_{k}=\beta\tau_{k}$
and $B_{k}=\overline{B}_{t_{k}}(x_{k})$. Clearly, the sequence
$\tau_{k}$ (and thus also $t_{k}$) does not increase. The balls
$\overline{B}_{\tau_{k}}(x_{k})$ are pairwise disjoint. Indeed, if
$l>k$ then $x_{l}\not\in B_{k}$, i.e., the pseudodistance between
$x_{l}$ and $x_{k}$ is greater than
$\beta\tau_{k}>2\tau_{k}\geq\tau_{k}+\tau_{l}$. Thus
$\overline{B}_{\tau_{k}}(x_{k})\cap\overline{B}_{\tau_{l}}(x_{l})=\emptyset$
by the triangle inequality for $d$. Now,
$$
\sum_{k\geq 1}\phi(\gamma t_{k})<\sum_{k\geq
1}\phi(\alpha\tau_{k})\leq\sum_{k\geq
1}\phi(\tau(x_{k}))=\sum_{k\geq
1}\xi(\overline{B}_{\tau_{k}}(x_{k}))\leq A;
$$
consequently, $\tau_{k}\to 0$, i.e., for each point $x$, not
belonging to the union of the balls $B_{k}$, $\tau(x)=0$, i.e., $x$
is a regular point. In addition, $t_{k}=\beta\tau_{k}\to 0$.
\end{proof}
\begin{remark}\label{re1}
{\rm (1) According to the Caratheodory construction, see, e.g., [F,
Ch.\ $\!$2.10], there is a finite Borel measure on $X$ whose
restriction to $\mathcal{F}$ is $\xi$.\\
(2) Assume that $\xi$ is the restriction to $\mathcal{F}$ of a Borel
measure $\mu$ on $X$ with support $\{x_{1},\dots, x_{n}\}$. Then, as
it follows from the proof, the number of the balls $B_{k}$ in this
case is $\leq n$ and the balls $\overline{B}_{\tau_{k}}(x_{k})$,
$k\geq 1$, cover the support of $\mu$. For otherwise, there is
$\overline{B}_{\tau_{k}}(x_{k})$ which does not meet $\{x_{1},\dots,
x_{n}\}$. Then $\overline{B}_{\tau(x_{k})}(x_{k})$ does not meet
$\{x_{1},\dots, x_{n}\}$, as well. Consequently,
$\mu(\overline{B}_{\tau(x_{k})}(x_{k}))=0$, a contradiction with the
choice of $x_{k}$.}
\end{remark}

Let $X$ be a pseudometric space with pseudometric $d$.  For every
$x\in X$ we set $S_{x}:=\{y\in X\ :\ d(x,y)=0\}$. Let $\mu$ be a
Borel measure on $X$ with $\mu(X)=k<\infty$ such that
$$
\int_{X}\ln^{+} d(x,\xi)\ \!d\mu(\xi)<\infty\ \ \ {\rm for\
all}\ \ \ x\in X,
$$
where $\ln^{+}t:=\max(0,\ln t)$.  Then we define
$$
u(x)=\left\{
\begin{array}{ccc}
\displaystyle \int_{X}\ln d(x,\xi)\ \!d\mu(\xi),&{\rm
if}&\mu(S_{x})=0\\
\\
\displaystyle -\infty,&{\rm if}&\mu(S_{x})>0.
\end{array}
\right.
$$
By definition, every Lebesgue integral $\int_{X}\ln d(x,\xi)\
\!d\mu(\xi)$ exists but may be equal to $-\infty$. In this case we
define $u(x)=-\infty$.
\begin{corollary}\label{cor1}
Fix $\gamma\in (0,1/2)$. Given $H>0$, $s>0$ there is a family of
closed balls $B_{j}$ with radii $r_{j}$ satisfying
$$
\sum r_{j}^{s}<\frac{(H/\gamma)^{s}}{s}
$$
such that
$$
u(x)\geq k\ln\left(\frac{H}{e}\right)\ \ \ {\rm for\ all}\ \ \ x\in
X\setminus\bigcup_{j}B_{j}.
$$
\end{corollary}
\begin{proof} Let $\phi(t)=(pt)^{s}$ be a majorant with
$p=\frac{(ks)^{1/s}}{H}$. We cover all irregular points of $X$ by
closed balls according to Lemma \ref{gorin} (with $\xi=\mu$) and
prove that the required inequality is valid for any regular point
$x$. This will complete the proof. First, observe that $u(x)$ is
finite for every regular point $x$ by the definition of the Lebesgue
integral and the regularity condition for the $\phi$. Let
$n(t;x)=\mu(\overline{B}_{t}(x))$ for such $x$. Then, for any
$N\geq\max(1,H)$ we have
$$
u(x)\geq\int_{\overline{B}(x,\xi)}\ln d(x,\xi)\
\!d\mu(\xi)=\int_{0}^{N}\ln t\ \!dn(t;x)=n(t;x)\ln
t|_{0}^{H}-\int_{0}^{N}\frac{n(t;x)}{t}\ \!dt.
$$
Since $n(t;x)<(pt)^{s}$, we obtain
$$
u(x)\geq n(N;x)\ln N-\int_{0}^{N}\frac{n(t;x)}{t}\ \!dt.
$$
In addition, $n(t;x)\leq n(N;x)$ for $t\leq N$. Therefore,
$$
\begin{array}{c}
\displaystyle u(x)\geq n(N;x)\ln N-\int_{0}^{H}\frac{(pt)^{s}}{t}\
\!dt-\int_{H}^{N}\frac{n(N;x)}{t}\ \!dt=\\
\\
\displaystyle n(N;x)\ln N-\frac{(pH)^{s}}{s}-n(N;x)\ln N+n(N;x)\ln H=
-k+n(N;x)\ln H.
\end{array}
$$
Letting here $N\to\infty$ and taking into account that
$\lim_{N\to\infty}n(N;x)=k$ we obtain the required result.
\end{proof}

We use also the following result proved in [L].
\begin{corollary}\label{cor2}
Let $f$ be a holomorphic function in the disk $|z|\leq 2eR$ $(R>0)$
in $\Co$, $f(0)=1$ and $\eta$ is an arbitrary positive number
$\leq\frac{3}{2}e$. Then inside the disk $|z|\leq R$ but outside a
family of closed disks $\overline{D}_{r_{i}}(z_{i})$ centered at
$z_{i}$ of radii $r_{i}$ such that $\sum r_{i}\leq 4\eta R$,
$$
\ln |f(z)|\geq -H(\eta)\ln M(2eR)
$$
where
$$
H(\eta)=2+\ln\left(\frac{3e}{2\eta}\right)
$$
and
$$
M(2eR):=\sup_{|z|\leq 2eR}|f|.
$$
\end{corollary}
\begin{remark}{\rm The proof is based on a particular case of Corollary
\ref{cor1} for $\mu$ a sum of delta-measures, and the Harnack
inequality for positive harmonic functions. According to Remark
\ref{re1} (2), from the proof presented in [L] it follows that the
number of disks $\overline{D}_{r_{i}}(z_{i})$ does not exceed the
number of zeros of $f$ in the disk $|z|<2R$ (which, by the Jensen
inequality, is bounded from above by $[\ln M_{f}(2eR)]$) and,
moreover, the disks $\overline{D}_{r_{i}/2}(z_{i})$ cover the set of
these zeros.}
\end{remark}
\subsubsection{}
In this part we will prove Theorem \ref{remez1}.
\begin{proof}
Let $X\subset\Co^{n}$ be a closed subset of the class
$\mathcal{A}_{2n}(s,a)$ where $s=2n-2+\alpha$, $\alpha>0$. Let $p$
be a holomorphic polynomial on $\Co^{n}$ of degree $k$. Let
$B\subset\Co^{n}$ be a closed Euclidean ball and $\omega\subset
X\cap B$ be an $\mathcal{H}_{s}$-measurable subset. We must prove the
inequality
\begin{equation}\label{ine0}
\sup_{B}|p|\leq\left(\frac{c_{1}\mathcal{H}_{2n}(B)}{\{\mathcal
{H}_{s}(\omega)\}^{2n/s}}\right)^{c_{2}k}\sup_{\omega}|p|
\end{equation}
where $c_{1}$ depends on $a$, $n$, $k$, $\alpha$ and $c_{2}>0$
depends on $\alpha$.

Since the ratio on the right-hand side of (\ref{ine0}) is invariant
with respect to dilations and translations of $\Co^{n}$ and the
class $\mathcal{A}_{2n}(s,a)$ is also invariant with respect to
these transformations, without loss of generality we may assume that
$B$ is the closed unit ball centered at $0\in\Co^{n}$. Then we must
prove that
\begin{equation}\label{ine1}
\sup_{B}|p|\leq\left(\frac{\overline
c_{1}}{\lambda^{2n/s}}\right)^{c_{2}k}\sup_{\omega}|p|
\end{equation}
where $\lambda:=\mathcal{H}_{s}(\omega)$, $\overline c_{1}$ depends
on $a$, $n$, $k$, $\alpha$ and $c_{2}>0$ depends on $\alpha$.

By $Z_{p}\subset\Co^{n}$ we denote the set of zeros of $p$.
According to Theorem \ref{dset}, $Z_{p}\in \mathcal{
A}_{2n}(2n-2,a,b)$ for some $a$ and $b$ depending on $n$ and $k$
only. By $\mathcal{H}_{2n-2,\ \!p}$ we denote the Hausdorff
$(2n-2)$-measure supported on $Z_{p}$. Let $B_{1}\subset B_{2}$ be
closed Euclidean balls centered at $0\in\Co^{n}$ of radii 2 and 10,
respectively. Set
$$
\mu:=\mathcal{H}_{2n-2,\ \!p}|_{B_{2}}.
$$
Since $Z_{p}\in \mathcal{A}_{2n}(2n-2,a,b)$, we have
\begin{equation}\label{con1}
\mu(B_{r}(x))\geq br^{2n-2}\ \ \ {\rm for\ all}\ \ \ x\in B_{1},\ \
\ 0\leq r\leq 5.
\end{equation}

Let $H>0$. Consider $\phi(t):=\frac{t^{s}}{H}$ as the majorant in
Lemma \ref{gorin}. Then a point $x\in\mathbf{C}^{n}$ is regular with
respect to $\phi$ and $\mu$ if $\mu(B_{r}(x))<\frac{r^{s}}{H}$ for
all $r>0$. (Here we consider $\Co^{n}$ with the Euclidean norm
$|\cdot|$.)
\begin{lemma}\label{lem1}
There is a sequence of open Euclidean balls $B_{r_{k}}(x_{k})$,
$k=1,2,\dots$, which collectively cover all the irregular points
such that
$$
\sum_{k\geq 1}r_{k}^{s}<3H\mu(B_{2}).
$$
Moreover, the distance $d(x)$ from a regular point $x$ to the
compact set $K:=B_{1}\cap Z_{p}$ is $\geq\min
\left\{5,\left(\frac{bH}{2^{s}}\right)^{1/\alpha}\right\}$.
\end{lemma}
\begin{proof} The first statement follows directly from Lemma \ref{gorin}.
Let $y\in K$ be such that $|x-y|=d(x)$. Observe that condition
(\ref{con1}) implies that $x\not\in K$. For otherwise, we must have
$$
br^{2n-2}<\frac{r^{s}}{H}\ \ \ {\rm for\ all}\ \ \ 0<r<5
$$
which is impossible. Thus $d(x)>0$. Next, the ball centered at $x$
of radius $2d(x)$ contains the ball centered at $y$ of radius
$d(x)$. Now from the regularity condition for $x$ by (\ref{con1}) we
get
$$
b\min\{5,d(x)\}^{2n-2}\leq\mu(B_{2d(x)}(x))\leq\frac{
\{2d(x)\}^{2n-2+\alpha}}{H}.
$$
This implies that
$$
d(x)\geq\min\left\{5,\left(\frac{bH}{2^{s}}\right)^{1/\alpha}\right\}.
$$
\end{proof}

Continuing the proof of the theorem observe that by the definition
of $X$,
\begin{equation}\label{con2}
\lambda:=\mathcal{H}_{s}(\omega)\leq a2^{s}
\end{equation}
(because if $\omega\subset X\cap B\neq\emptyset$, then $\omega$ is
contained in a closed Euclidean ball of radius $2$ centered at a
point of $X$). Without loss of generality we may assume that
$\lambda>0$.
\begin{lemma}\label{lem2}
The set $\omega$ cannot be covered by a family $\{B_{j}\}$of open
Euclidean balls whose radii $r_{j}$ satisfy
$$
\sum r_{j}^{s}<\frac{\lambda}{2^{s}a}.
$$
\end{lemma}
\begin{proof} Assume to the contrary that there is a family of balls
$B_{j}:=B_{r_{j}}(x_{j})$, $j=1,2,\dots$, whose radii satisfy the
inequality of the lemma which covers $\omega$. Without loss of
generality we may assume that each $B_{j}$ meets $\omega$. Then for
every $x_{j}$ choose $y_{j}\in \omega$ so that $|x_{j}-y_{j}|\leq
r_{j}$. Clearly, the family of balls $\{B_{2r_{j}}(y_{j})\}$ also
covers $\omega$. From here, since $\omega\subset X\in \mathcal{
A}_{2n}(s,a)$, we obtain
$$
\lambda:=\mathcal{H}_{s}(\omega)\leq\sum\mathcal{H}_{s}(X\cap
B_{2r_{j}}(y_{j}))\leq 2^{s}a\sum r_{j}^{s}<\lambda,
$$
a contradiction.
\end{proof}

Further, note that $\mu(B_{2})$ in Lemma \ref{lem1} is bounded from
above by a constant $c$ depending on $n$ and $k$ only (because
$Z_{p}\in \mathcal{A}_{2n}(2n-2,a,b)$ with $a$, $b$ depending on
$n$, $k$ only). Thus choosing in this lemma $H:=\frac{\lambda}{3c
2^{s}}$ we obtain from Lemma \ref{lem2} for some $\overline c$
depending on $n$, $k$:
\begin{corollary}\label{c3}
There is a point $x\in\omega$ such that
$$
dist(x,Z_{p})\geq\min\left\{1,(\overline{c}\lambda)^{1/\alpha}\right\}.
$$
\end{corollary}
\begin{proof} From the above lemmas it follows that there is $x\in\omega$
such that
$$
dist(x,Z_{p}\cap
B_{1})\geq\min\left\{5,(\overline{c}\lambda)^{1/\alpha}\right\}.
$$
Moreover, $x\in B$ and so $dist(x,Z_{p}\setminus B_{1})\geq 1$; this
implies the required.
\end{proof}

Let $z\in B$ be a point such that
$$
M:=\max_{B}|p|=|p(z)|.
$$
Let $l$ be the complex line passing through $z$ and the point $x$
from Corollary \ref{c3}. Without loss of generality we may identify
$l$ with $\Co$ so that $z$ coincides with $0\in\Co$. Then, in this
identification, the point $x$ belongs to $\overline D_{2}(0)$, the
closed disk of radius 2 centered at $0$. Observe also that (under
the identification) the set $B_{1}\cap l$ contains
$\overline{D}_{1}(0)$. Thus, by the classical Bernstein inequality
for holomorphic polynomials
$$
\max_{|z|\leq 4e}|p|\leq (4e)^{k}\max_{|z|\leq 1}|p|\leq
(4e)^{k}\max_{B_{1}}|p|\leq (8e)^{k}\max_{B}|p|:=(8e)^{k}M.
$$

Set $f=p/M$ and apply Corollary \ref{cor2} with $R=2$. According to
this corollary for every $\eta\leq 3e/2$ there is a family of closed
disks $\overline D_{r_{i}}(z_{i})$ such that $\sum r_{i}\leq 8\eta$
and $\ln |f(z)|>-H(\eta)k\ln (8e)$ for any $|z|\leq 2$
outside the above disks where $H(\eta)=2+\ln (3e/2\eta)$.
Recall also that the number of these disks is $\leq$ the number of
zeros of $f$ in $|z|\leq 4$ and the disks $\overline
D_{r_{i}/2}(z_{i})$ cover the set of zeros of $f$ there. In
particular, if a point $z\in \overline D_{1}(0)$ satisfies
$dist(z,Z_{f})\geq 14\eta$ where $Z_{f}$ is the set of zeros of $f$
in $\Co$, then it cannot belong to the union of disks $\overline
D_{r_{i}}(z_{i})$, and therefore $|f(z)|$ satisfies the above
inequality. Choose $\eta:=\min(1, (\overline
c\lambda)^{1/\alpha})/14$. Then by Corollary \ref{c3},
$dist(x,Z_{f})\geq 14\eta$. Thus we have
$$
\sup_{\omega}\ln |f|\geq\ln |f(x)|\geq -H(\eta)k\ln (8e).
$$
We will consider two cases:

(1)
$$
(\overline c\lambda)^{1/\alpha}\geq 1. \
$$
Then $\eta=\frac{1}{14}$ and
$$
\sup_{\omega}\ln |f|\geq - (3+\ln 21)k\ln (8e)>-20d.
$$
This and (\ref{con2}) imply that
$$
\sup_{B}|p|\leq
e^{20k}\sup_{\omega}|p|=\left(\frac{e^{20}}{\lambda^{2n/s}}\lambda^{2n/s}\right)^{k}
\sup_{\omega}|p|\leq\left(\frac{2^{2n}a^{2n/s}e^{20}}{\lambda^{2n/s}}\right)^{k}
\sup_{\omega}|p|.
$$
Thus, inequality (\ref{ine1}) is proved in this case.

(2)
$$
(\overline c\lambda)^{1/\alpha}< 1.
$$
Then
$$
\sup_{\omega}\ln |f|\geq -(c'-\ln \lambda^{1/\alpha})k\ln (8e)
$$
where $c'$ depends on $n$ and $k$ only. This yields
$$
\sup_{B}|p|\leq\left(\frac{\overline
c_{1}}{\lambda^{2n/s}}\right)^{c_{2}k}\sup_{\omega}|p|
$$
where $\overline c_{1}>0$ depends on $k$, $n$ and $\alpha$ and
$c_{2}>0$ depends on $\alpha$ only.

The proof of Theorem \ref{remez1} is complete.
\end{proof}
\subsubsection{Proof of Corollary \ref{remez2}.} The proof follows
directly from the estimates obtained in cases (1) and (2) above  and
from the
fact that $X\in \mathcal{A}_{2n}(s,a,b)$.\ \ \ \ \ $\Box$\\
\subsubsection{Proofs of Corollaries \ref{bmo} and \ref{holder}.}
The proofs repeat word-for-word proofs of similar statements of
Theorems 1 and 3 of [ABr] and are based on the inequality of
Corollary \ref{remez2}. We leave the details to the readers.\ \ \ \
\ $\Box$
\section{Real Polynomials}
\subsection{Weak Remez Type Inequalities} In this section we will prove
Theorem \ref{weak}.
\begin{proof}
We set for brevity
$$
||p;\omega||_{q}=\left(\frac{1}{\mathcal{H}_{s}(\omega)}\int_{\omega}|p|^{q}\
\!d\mathcal{H}_{s} \right)^{1/q}\ \ \  {\rm and}
$$
$$
||p;U||_{r}=\left(\frac{1}{\mathcal{H}_{n}(U)}\int_{U}|p|^{r}\
\!d\mathcal{H}_{n} \right)^{1/r}.
$$

Since the above functions are invariant with respect to dilations of
$\Re^{n}$, without loss of generality we may and will assume that
$\mathcal{H}_{n}(U)=1$.

Let $\Sigma(a,\lambda)$, $a,\lambda>0$, be the class of subsets
$\omega\in \mathcal{A}_{n}(s,a)$ of $U$ satisfying
\begin{equation}\label{ratio1}
\{\mathcal{H}_{s}(\omega)\}^{n/s}\geq\lambda.
\end{equation}
We must show that there is a positive constant
$C=C(U,n,q,r,s,k,a,\lambda)$ such that for every real polynomial $p$
of degree $k$ on $\Re^{n}$
\begin{equation}\label{weak2}
||p;U||_{r}\leq C||p;\omega||_{q}.
\end{equation}
\begin{remark}
{\rm Let $C_{0}$ be the optimal constant in (\ref{weak2}). Since the
class $\Sigma(a,\lambda)$ increases as $\lambda$ decreases, $C_{0}$
increases in $1/\lambda$, as is required in the theorem.}
\end{remark}

If, on the contrary, the constant in (\ref{weak2}) does not exist,
one can find sequences of real polynomials $\{p_{j}\}$ of degrees
$k$ and sets $\{\omega_{j}\}\subset\Sigma(a,\lambda)$ so that
\begin{equation}\label{weak3}
||p_{j};U||_{r}=1\ \ \ {\rm for\ all}\ \ \ j\in\N,
\end{equation}
\begin{equation}\label{weak4}
\lim_{j\to\infty}||p_{j};\omega_{j}||_{q}=0.
\end{equation}
Since all (quasi-) norms on the space of real polynomials of degree
$k$ on $\Re^{n}$ are equivalent, (\ref{weak3}) implies the existence
of a subsequence of $\{p_{j}\}$ that converges uniformly on $U$ to a
polynomial $p\in\Re[x_{1},\dots, x_{n}]$ with $deg\ \!p\leq k$.
Assume without loss of generality that $\{p_{j}\}$ itself converges
uniformly to $p$. Then (\ref{weak3}), (\ref{weak4}) imply for this
$p$ that
\begin{equation}\label{weak5}
||p;U||_{r}=1,
\end{equation}
\begin{equation}\label{weak6}
\lim_{j\to\infty}||p;\omega_{j}||_{q}=0.
\end{equation}
From this we derive the next result.
\begin{lemma}\label{lweak}
There is a sequence of closed subsets
$\{\sigma_{j}\}\subset\overline{U}$ such that for every $j$ larger
than a fixed $j_{0}$ the following is true
\begin{equation}\label{weak7}
\{\mathcal{H}_{s}(\sigma_{j})\}^{n/s}\geq\frac{\lambda}{2^{n/s}}.
\end{equation}
Moreover,
\begin{equation}\label{weak8}
\max_{\sigma_{j}}|p|\to 0\ \ \ {\rm as}\ \ \ j\to\infty.
\end{equation}
\end{lemma}
\begin{proof}
Let first $q<\infty$. By the (probabilistic) Chebyshev inequality
$$
\mathcal{H}_{s}(\{x\in\omega_{j}\ :\ |p(x)|\leq
t\})\geq\mathcal{H}_{s}(\omega_{j})-\frac{\mathcal{H}_{s}(\omega_{j})}{t^{q}}
||p;\omega_{j}||_{q}^{q}.
$$
Pick here $t=t_{j}:=||p;\omega_{j}||_{q}^{1/2}$. Then by
(\ref{weak6}) the left-hand side is at least
$\frac{1}{2}\mathcal{H}_{s}(\omega_{j})$, for $j$ sufficiently
large. Denoting the closure of the set in the braces by $\sigma_{j}$
we also have
$$
\max_{\sigma_{j}}|p|=t_{j}\to 0\ \ \ {\rm as}\ \ \ j\to\infty.
$$
If $q=\infty$, simply set $\sigma_{j}:=\omega_{j}$ to produce
(\ref{weak7}) and (\ref{weak8}).
\end{proof}

Apply now the Hausdorff compactness theorem to find a subsequence of
$\{\sigma_{j}\}$ converging to a closed subset
$\sigma\subset\overline{U}$ in the Hausdorff metric. We assume
without loss of generality that $\{\sigma_{j}\}\to\sigma$. By
(\ref{weak8}) this limit set is a subset of the zero set of $p$.
Since $p$ is nontrivial by (\ref{weak5}), the dimension of its zero
set is at most $n-1$; hence $\mathcal{H}_{s}(\sigma)=0$ because
$s>n-1$. Then for every $\epsilon>0$ one can find a finite open
cover of $\sigma$ by open Euclidean balls $B_{i}$ of radii $r_{i}$
at most $r(\epsilon)$ so that
\begin{equation}\label{weak9}
\sum_{i}r_{i}^{s}<\epsilon.
\end{equation}
Let $\sigma_{\delta}$ be a $\delta$-neighbourhood of $\sigma$ such
that
$$
\sigma_{\delta}\subset\bigcup_{i}B_{i}\ \ \ {\rm and}\ \ \
\delta<r(\epsilon).
$$
Pick $j$ so large that $\sigma_{j}\subset\sigma_{\delta}$. For every
$B_{i}$ intersecting $\sigma_{j}$ choose a point $x_{i}\in
B_{i}\cap\sigma_{j}$. Consider an open Euclidean ball $\widetilde
B_{i}$ centered at $x_{i}$ of radius twice that of $B_{i}$. Then
$B_{i}\subset\widetilde B_{i}$ and $\{\widetilde B_{i}\}$ is an open
cover of $\sigma_{j}$. Hence
$$
\mathcal{H}_{s}(\sigma_{j})\leq\sum_{i}\mathcal{H}_{s}(\sigma_{j}\cap\widetilde
B_{i})\leq a2^{s}\sum_{i}r_{i}^{s}
$$
because $\omega_{j}\in\mathcal{A}_{n}(s,a)$. Together with
(\ref{weak7}) and (\ref{weak9}) this implies that
$$
\frac{1}{2}\lambda^{s/n}\leq\mathcal{H}_{s}(\sigma_{j})\leq a
2^{s}\sum_{i}r_{i}^{s}\leq 2^{s}a\epsilon.
$$
Letting $\epsilon\to\infty$ one gets a contradiction.
\end{proof}
\subsection{Strong Remez Type Inequalities}
Strong Remez type inequalities for real polynomials from $\Re[x]$ and Ahlfors regular subsets of $\Re$
are proved in [ABr]. Inequalities of the form described in Theorem
\ref{remez1} are also valid for real polynomials on $\Re^{2}$. The
method of the proof of such inequalities is very similar to that of
Theorem \ref{remez1} and is based on the fact that an analytic
compact connected curve in $\Re^{n}$ is a $1$-set. (The detailed
proof will be presented elsewhere.) It is still an open question
whether similar strong Remez type inequalities are valid for real
polynomials on $\Re^{n}$ for $n>2$.
\section{Proof of Theorem \ref{ybr}}
\begin{proof}
It is well known, see, e.g., [JW, Prop.$\ \!$VIII.1], that the closure
$\overline X$ of an $s$-set $X$ is also an $s$-set and $\mathcal{H}_{s}(\overline X\setminus X)=0$. Moreover, the spaces 
$\dot{C}_{q}^{k,\omega}(\overline X)$ 
and $\dot{C}_{q}^{k,\omega}(X)$ are isometric. Thus without loss of generality
we may and will assume in the proof that $X$ is closed. 

Given $f\in\dot{C}_{q}^{k,\omega}(X)$ we should find a function
$\widetilde f:X\to\Re$ which equals $f$ modulo zero 
$\mathcal{H}_{s}$-measure and
admits an extension to a function from $\dot{\Lambda}^{k,\omega}(\Re^{n})$.

We begin with
\begin{lemma}\label{l5.1}
Let $\omega:\Re_{+}\to\Re_{+}$ be a quasipower $k$-majorant, see
Definition \ref{morrey5}. Let $t_{j}:=2^{j}$, $j\in\Z_{+}$. Then for
every pair of integers $-\infty<i<i'<\infty$ we have
\begin{equation}\label{e5.2}
\sum_{j=i}^{i'}\omega(t_{j})\leq c(k,\omega)\ \!\omega(t_{i'}).
\end{equation}
\end{lemma}
\begin{proof}
By the monotonicity of $\omega$
$$
\omega(t_{j})\leq\frac{1}{\ln 2}\int_{t_{j}}^{t_{j+1}}
\frac{\omega(u)}{u}\ \!du
$$
and therefore the sum in (\ref{e5.2}) is at most
$$
\frac{1}{\ln 2}\int_{t_{i}}^{t_{i'+1}}\frac{\omega(u)}{u}\ \!du\leq
\frac{1}{\ln 2}\ \!C_{\omega}\ \!\omega(t_{i'+1})
\leq\frac{C_{\omega}}{\ln 2}2^{k}
\omega(t_{i'}).
$$
\end{proof}

Our next result reformulates a theorem of the paper [BSh1], see also
[BSh2, $\!$Th. $\!$3.5] concerning the trace of the space 
$\dot{\Lambda}^{k,\omega}(\Re^{n})$ to an
arbitrary closed subset $X\subset\Re^{n}$, to adopt it to our situation.
The trace space denoted by
$\dot{\Lambda}^{k,\omega}(\Re^{n})|_{X}$ consists of locally bounded
functions $f:X\to\Re$ and is equipped with seminorm
\begin{equation}\label{e5.6}
|f|_{\dot{\Lambda}^{k,\omega}(\Re^{n})|_{X}}:=
\inf\{|g|_{\dot{\Lambda}^{k,\omega}(\Re^{n})}\ :\ f=g|_{X}\}.
\end{equation}

To formulate the result we need
\begin{definition}\label{d5.2}
Let $X\subset\Re^{n}$ and $\omega:\Re_{+}\to\Re_{+}$ be as above, and
$\mathcal{T}_{\omega}:=\{t_{i}\}_{i\in\Z_{+}}$ be the sequence of Lemma
\ref{l5.1}.

A family $\Pi:=\{P_{Q}\}_{Q\in\mathcal{K}_{X}}$ of polynomials of degree
$k-1$ is said to be a $(k,\omega,X)$-chain if for every pair of cubes
$Q\subset Q'$ from $\mathcal{K}_{X}$ which satisfy for some $i\in\Z$
the condition
\begin{equation}\label{e5.7}
t_{i}\leq r_{Q}<r_{Q'}\leq t_{i+2}
\end{equation}
the inequality
\begin{equation}\label{e5.8}
\max_{x\in Q}|P_{Q}(x)-P_{Q'}(x)|\leq C\omega(r_{Q'})
\end{equation}
holds with a constant $C$ independent of $Q$, $Q'$ and $i$.
\end{definition}

The linear space of such chains is denoted by $Ch(k,\omega,X)$. It is
equipped with seminorm
$$
|\Pi|_{Ch}:=\inf C
$$
where the infimum is taken over all constants $C$ in (\ref{e5.8}).

Recall that $\mathcal{K}_{X}$ is the family of closed cubes centered at $X$ and of radii at most $4 diam\ \!X$. 
In the sequel $c_{Q}$ and $r_{Q}$ stand for
the center and the radius of the cube $Q$.

Using the concept introduced and the related notations we now formulate
the desired result.
\begin{proposition}\label{p5.3}
(a) A locally bounded function $f:X\to\Re$ belongs to 
$\dot{\Lambda}^{k,\omega}(\Re^{n})|_{X}$ if and only if there is a
$(k,\omega,X)$-chain $\Pi:=\{P_{Q}\}_{Q\in\mathcal{K}_{X}}$ such that
for every $Q\in\mathcal{K}_{X}$
\begin{equation}\label{e5.9}
f(c_{Q})=P_{Q}(c_{Q}).
\end{equation}
Moreover, the following two-sided inequality
$$
|\Pi|_{Ch}\approx |f|_{\dot{\Lambda}^{k,\omega}(\Re^{n})|_{X}}
$$
holds with constants independent of $f$.\\
(b) If, in addition, this chain depends on $f$ linearly, then there is
a linear extension operator $T_{k}:\dot{\Lambda}^{k,\omega}(\Re^{n})|_{X}
\to\dot{\Lambda}^{k,\omega}(\Re^{n})$ such that
$$
||T_{k}||\leq O(1)|\Pi|_{Ch}.
$$
\end{proposition}
Hereafter $O(1)$ denotes a constant depending only on inessential
parameters. It may change from line to line and even in a single line.
\begin{proof}
In the above cited papers this result is proved under the assumption
that inequality
(\ref{e5.8}) holds for any pair of cubes $Q\subset Q'$ centered at
$X$. The restrictions (\ref{e5.7}) and 
$r_{Q},\ r_{Q'}\leq 4 diam\ \!X$ may be not satisfied for this pair.
In the forthcoming derivation we explain how these restrictions can
be disregarded to apply the aforementioned Theorem 3.5 of [BSH2] and in
this way to complete the proof of the proposition.

Consider first the case of an unbounded $X$. Hence, the only restriction
is now inequality (\ref{e5.7}) and we should show that
if a $(k,\omega,X)$-chain satisfies
condition (\ref{e5.8}) under restriction (\ref{e5.7}), then
(\ref{e5.8}) holds for {\em any} pair $Q\subset Q'$ from $\mathcal{K}_{X}$.
Note that the necessity of conditions (\ref{e5.8}) and (\ref{e5.9})
trivially follows from that in the aforementioned Theorem 3.5 from [BSh2]. So
we should only prove their sufficiency.

Assume that $f\in l_{\infty}^{loc}(X)$ and conditions 
(\ref{e5.7})-(\ref{e5.9}) hold. Let $Q\subset Q'$ be a pair of cubes from 
$\mathcal{K}_{X}$ of radii $r$ and $r'$, respectively. Then for some
indices $i\leq i'$
$$
t_{i}\leq r\leq t_{i+1}\ \ \ {\rm and}\ \ \ t_{i'}\leq r'\leq t_{i'+1}.
$$
If $i=i'$, then by (\ref{e5.8})
$$
\max_{Q}|P_{Q}-P_{Q'}|\leq 2|\Pi|_{Ch}\omega(t_{i+1})\leq
2\left(\frac{t_{i+1}}{t_{i}}\right)^{k}\omega(r')|\Pi|_{Ch}=
2^{k+1}\omega(r')|\Pi|_{Ch}
$$
as is required.

Let now $i<i'$ and $r_{j}$ with $i\leq j\leq i'+1$ are given by
$$
r_{i}:=r,\ \ \ r_{i'+1}=2r'\ \ \ {\rm and}\ \ \ r_{j}:=t_{j}\ \ \
{\rm for}\ \ \ i<j<i'+1.
$$
Let $Q_{j}$ be the cubes centered at $c_{Q}$ of radii
$r_{j}$, $i\leq j<i'+1$, and $Q_{i'+1}$ be the cube centered at $c_{Q'}$ of radius $r_{i'+1}$. (In particular, $\{Q_{j}\}_{i\leq j\leq i'+1}\subset\mathcal{K}_{X}$ is an increasing sequence of cubes with $Q_{i}:=Q$.) Then 
\begin{equation}\label{e5.10}
\max_{Q}|P_{Q}-P_{Q'}|\leq\sum_{j=i}^{i'}\ \max_{Q_{j+1}}|P_{Q_{j}}-
P_{Q_{j+1}}|.
\end{equation}
It is easily seen that (\ref{e5.7})
holds for every pair $Q_{j}\subset Q_{j+1}$, $i\leq j\leq i'$. Applying 
(\ref{e5.8}) to each of these pairs and then (\ref{e5.2}) and the definition of $\omega$ 
we estimate the right-hand side of (\ref{e5.10}) by
$$
2|\Pi|_{Ch}\sum_{j=i}^{i'}\omega(r_{j+1})\leq O(1)|\Pi|_{Ch}\omega(t_{i'+2})
\leq O(1)|\Pi|_{Ch}\omega(r').
$$

Thus we conclude that inequality (\ref{e5.8}) holds for every pair
$Q\subset Q'$ of cubes centered at $X$.

Let now $diam\ \!X<\infty$. The previous argument proves the required
inequality
\begin{equation}\label{nabla}
\max_{Q}|P_{Q}-P_{Q'}|\leq C\omega(r_{Q'})
\end{equation}
for every pair $Q\subset Q'$ from $\mathcal{K}_{X}$ under
the restriction $r_{Q'}\leq 2diam\ \!X$.
Fix a cube $\widetilde Q\in\mathcal{K}_{X}$ with $r_{\widetilde Q}=2diam\ \!X$ and introduce a new family of polynomials
$\{\overline P_{Q}\}$, where $Q$ runs over the set of all cubes centered at $X$, by
setting
\begin{equation}\label{fnew}
\overline P_{Q}:=\left\{
\begin{array}{ccc}
P_{Q},&{\rm if}&r_{Q}\leq diam\ \!X\\
\\
P_{\widetilde Q}-P_{\widetilde Q}(c_{Q})+f(c_{Q}),&{\rm if}& r_{Q}>diam\ \!X.
\end{array}
\right.
\end{equation}
We will prove that the new family satisfies the hypotheses of Theorem 3.5 from [BSh2]. This will complete the proof of the proposition in this case.

Clearly, $\{\overline P_{Q}\}$ satisfies condition (\ref{e5.9}), and
if the chain $\Pi$ depends linearly on $f$, then $\{\overline P_{Q}\}_{Q}$ depends linearly on $f$, as well. So we must check only that (\ref{e5.8}) holds for $\{\overline P_{Q}\}$ for every pair $Q\subset Q'$ of cubes centered at $X$. According to (\ref{nabla}) and (\ref{fnew}) inequality (\ref{e5.8}) holds for this family for every pair of cubes $Q\subset Q'$
with $r_{Q'}\leq diam\ \!X$. Assume now that $r_{Q'}\geq r_{Q}>diam\ \!X$.
Then by (\ref{fnew}) we have
$$
\begin{array}{c}
\displaystyle
\max_{Q}|\overline P_{Q}-\overline P_{Q'}|\leq |P_{\widetilde Q}(c_{Q})-f(c_{Q})|+
|P_{\widetilde Q}(c_{Q'})-f(c_{Q'})|\leq\\
\\
\displaystyle
\max_{Q_{1}}|P_{\widetilde Q}-P_{Q_{1}}|+\max_{Q_{2}}|P_{\widetilde Q}-P_{Q_{2}}|\leq 2C\omega(r_{\widetilde Q})\leq O(1)C\omega(r_{Q'}).
\end{array}
$$
Here $Q_{1}$ and $Q_{2}$ are some cubes from $\mathcal{K}_{X}$ centered at $c_{Q}$ and $c_{Q'}$, respectively, and contained in $Q$. The last two inequalities follow from (\ref{nabla}) and the definition of $\omega$. 
Finally, if $r_{Q}\leq diam\ \!X<r_{Q'}$, then $Q\subset\widetilde Q$ and so 
we have by (\ref{nabla}) and by the definition of $\omega$
$$
\max_{Q}|\overline P_{Q}-\overline P_{Q'}|\leq \max_{Q}|P_{Q}-P_{\widetilde Q}|+|P_{\widetilde Q}(c_{Q'})-f(c_{Q'})|\leq 2C\omega(r_{\widetilde Q})\leq O(1)C\omega(r_{Q'})
$$
as is required.

Hence, in both of these cases the assumptions of Theorem 3.5 from [BSh2] hold.
This completes the proof of the proposition.
\end{proof}

Now we outline the proof of Theorem \ref{ybr}. Given 
$f\in\dot{C}_{q}^{k,\omega}(X)$ where $X\subset\Re^{n}$ is a closed $s$-set,
$n-1<s\leq n$, we will define a new function
$\widetilde f:X\to\Re$ such that
\begin{equation}\label{e5.11}
\widetilde f(x)=f(x)\ \ \ \mathcal{H}_{s}-almost\ everywhere\ on\ X.
\end{equation}
We then apply Proposition \ref{p5.3} to this function
to show that $\widetilde f\in\dot{\Lambda}^{k,\omega}(\Re^{n})|_{X}$
to construct a linear extension operator from
$\dot{C}_{q}^{k,\omega}(X)$ to $\dot{\Lambda}^{k,\omega}(\Re^{n})$.
To this end we will find for the $\widetilde f$ a $(k,\omega,X)$-chain
linearly depending on $f$. 
In the definition of the desired chain we will use the following 
construction. Let 
$Q:=Q_{r}(x)\in\mathcal{K}_{X}$. By the Kadets-Snobar theorem [KS] there
is a linear projection $\pi_{Q}$ from 
$L_{1}(X_{r}(x);\mathcal{H}_{s})$ onto the subspace of polynomials of degree
$k-1$ restricted to $X_{r}(x):=Q_{r}(x)\cap X$ whose norm 
$||\pi_{Q}||_{1}\leq \sqrt{d_{k,n}}$ where $d_{k,n}$ is the dimension of
the space of polynomials of degree $k-1$ on $\Re^{n}$. Set
\begin{equation}\label{e5.12}
P_{Q}(f):=\pi_{Q}(f).
\end{equation}
Using the definitions of $\widetilde f$ and 
$\{P_{Q}(f)\}_{Q\in\mathcal{K}_{X}}$ we will show that the following is true.
\begin{claim}
There exists a $(k,\omega,X)$-chain 
$\widetilde\Pi(f):=\{\widetilde P_{Q}(f)\}_{Q\in\mathcal{K}_{X}}$ linearly
depending on $f$ and such that
\begin{equation}\label{e5.13}
|\widetilde\Pi(f)|_{Ch}\leq O(1)|f|_{\dot{C}_{q}^{k,\omega}(\Re^{n})|_{X}}.
\end{equation}
\end{claim}
\begin{claim}
For every $Q\in\mathcal{K}_{X}$
\begin{equation}\label{e5.14}
\widetilde f(c_{Q})=\widetilde P_{Q}(f)(c_{Q}).
\end{equation}
\end{claim}
Since the operator $f\mapsto \widetilde P_{Q}(f)$ 
is linear, these allow us to apply
Proposition \ref{p5.3} and to conclude that 
$\widetilde f\in\dot{\Lambda}^{k,\omega}(\Re^{n})|_{X}$, and there is a linear
extension operator $T_{k}:\dot{C}_{q}^{k,\omega}(X)\to
\dot{\Lambda}^{k,\omega}(\Re^{n})$ satisfying
$$
||T_{k}||\leq O(1)
$$
completing the first part of the proof of Theorem \ref{ybr}. The fact that
the restriction to $X$ of every $f\in\dot{\Lambda}^{k,\omega}(\Re^{n})$ belongs to $\dot{C}^{k,\omega}_{q}(X)$ follows easily from Proposition \ref{p5.3} and Definition \ref{morrey4}. This proves also the second assertion of the
theorem and completes its proof.

To realize this program we need several auxiliary results. The main tool
in their proofs is the weak Remez type inequality for $s$-sets, see 
Theorem \ref{weak} and (\ref{weak1}).
\begin{lemma}\label{l5.4}
For every $Q=Q_{r}(x)\in\mathcal{K}_{X}$
\begin{equation}\label{e5.15}
\left\{\frac{1}{\mathcal{H}_{s}(X_{r}(x))}\int_{X_{r}(x)}
|f-\pi_{Q}(f)|^{q}\ \!d\mathcal{H}_{s}\right\}^{1/q}\leq O(1)
\mathcal{E}_{k}(f;Q).
\end{equation}
\end{lemma}
\begin{proof}
Here and below for $Q=Q_{r}(x)\in\mathcal{K}_{X}$ by 
$P_{Q}$ we denote a polynomial
of degree $k-1$ satisfying
\begin{equation}\label{e5.16}
\left\{\frac{1}{\mathcal{H}_{s}(X_{r}(x))}\int_{X_{r}(x)}
|f-P_{Q}|^{q}\ \!d\mathcal{H}_{s}\right\}^{1/q}=\mathcal{E}_{k}(f;Q).
\end{equation}
Then 
$$
f-\pi_{Q}(f)=(f-P_{Q})+\pi_{Q}(f-P_{Q}),
$$
and applying the triangle inequality we estimate the left-hand side in 
(\ref{e5.15}) as is required but with the factor $(1+||\pi_{Q}||_{q})$ 
instead of $O(1)$. So it remains to show that $||\pi_{Q}||_{q}\leq O(1)$.
However, for $q=1$ this norm is bounded by $\sqrt{d_{k,n}}$ by the
definition. On the other hand, the weak Remez type inequality, see
(\ref{weak1}), and the fact that $X$ is an $s$-set,  imply that
$$
||\pi_{Q}(g)||_{1}\approx ||\pi_{Q}(g)||_{q}
$$
with the constants of equivalence independent of $g$ and $Q$. 
Thus by the H\"{o}lder inequality we have
$$
||\pi_{Q}(g)||_{q}\leq O(1)||\pi_{Q}(g)||_{1}\leq O(1)||\pi_{Q}||_{1}||g||_{1}
\leq O(1)||g||_{q}
$$
\end{proof}
\begin{lemma}\label{l5.5}
Let $Q=Q_{r}(x)\in\mathcal{K}_{X}$. Then there exists the limit
\begin{equation}\label{e5.17}
\widetilde f(x):=\lim_{Q\to x}P_{Q}(x)
\end{equation}
and, moreover,
\begin{equation}\label{e5.18}
|\widetilde f(x)-P_{Q}(x)|\leq O(1)\omega(r)|f|_{\dot{C}_{q}^{k,\omega}(X)}.
\end{equation}
\end{lemma} 
\begin{proof}
Let $i$ be defined by
\begin{equation}\label{e5.19'}
t_{i}<r\leq t_{i+1}
\end{equation}
and for $j\leq i$
$$
Q_{j}:=Q_{t_{j}}(x),\ \ \ P_{j}:=P_{Q_{j}}.
$$
Recall that $\{t_{j}\}$ is the sequence of Lemma \ref{l5.1}. We also
set $Q_{i+1}:=Q$ and $P_{i+1}:=P_{Q}$. Since $X$ is an $s$-set, the
weak Remez type inequality (\ref{weak1}) implies that
$$
|P_{j+1}(x)-P_{j}(x)|\leq O(1)|||P_{j+1}-P_{j};X_{j}|||
$$
where for simplicity we set
$$
|||g;X_{j}|||:=\left\{\frac{1}{\mathcal{H}_{s}(X_{j})}
\int_{X_{j}}|g|^{q}\ \!d\mathcal{H}_{s}\right\}^{1/q}
\ \ \ {\rm and}\ \ \ X_{j}:=Q_{j}\cap X.
$$
Adding and subtracting $f$ and remembering the definition of $P_{j}$,
see (\ref{e5.16}), we estimate the right-hand side of the last inequality
by
$$
O(1)\{\mathcal{E}_{k}(f;Q_{j})+|||f-P_{j+1};X_{j}|||\}.
$$
By definition, the first term is bounded by $\omega(t_{j})
|f|_{\dot{C}_{q}^{k,\omega}(X)}$ while the second one is at most
$$
\left(\frac{\mathcal{H}_{s}(X_{j+1})}{\mathcal{H}_{s}(X_{j})}\right)^{1/q}
\mathcal{E}_{k}(f;Q_{j+1})\leq\left(\frac{at_{j+1}^{s}}{bt_{j}^{s}}
\right)^{1/q}\omega(t_{j+1})|f|_{\dot{C}^{k,\omega}_{q}(X)},
$$
see Definition \ref{d2} of $s$-sets.\\
Since, in turn, $t_{j+1}/t_{j}\leq 2$, using the definition of $\omega$
we finally get
$$
|P_{j+1}(x)-P_{j}(x)|\leq O(1)\omega(t_{j})|f|_{\dot{C}^{k,\omega}_{q}(X)}.
$$
This, Lemma \ref{l5.1} and the choice of $i$, see (\ref{e5.19'}), yield
$$
\begin{array}{c}
\displaystyle
\sum_{j\leq i}|P_{j+1}(x)-P_{j}(x)|\leq O(1)|f|_{\dot{C}^{k,\omega}_{q}(X)}
\sum_{j\leq i}\omega(t_{j})
\leq O(1)|f|_{\dot{C}^{k,\omega}_{q}(X)}\omega(t_{i})\leq \\
\\
\displaystyle
O(1)\omega(r)|f|_{\dot{C}^{k,\omega}_{q}(X)}.
\end{array}
$$
This implies easily that the limit
$$
\widetilde f(x):=\lim_{Q\to x}P_{Q}(x)=P_{i+1}(x)+\sum_{j\leq i}
(P_{j}(x)-P_{j+1}(x))
$$
exists and, moreover,
$$
|\widetilde f(x)-P_{Q}(x)|\leq O(1)\omega(r)|f|_{\dot{C}^{k,\omega}_{q}(X)}.
$$
\end{proof}
\begin{lemma}\label{l5.5a}
The assertions of the previous lemma hold with the same $\widetilde f(x)$
for $P_{Q}(f)$ substituted for $P_{Q}$.
\end{lemma}
\begin{proof}
By (\ref{e5.12})
$$
P_{Q}-P_{Q}(f)=\pi_{Q}(P_{Q}-f)
$$
and then Lemma \ref{l5.4} and inequality (\ref{weak1}) yield
$$
\begin{array}{c}
\displaystyle
|P_{Q}(x)-P_{Q}(f)(x)|\leq O(1)\max_{Q\cap X}|P_{Q}-P_{Q}(f)|\leq
O(1)|||P_{Q}-P_{Q}(f);Q\cap X|||\leq
\\
\\
\displaystyle
O(1)\{\mathcal{E}_{k}(f;Q)+|||f-P_{Q};Q\cap X|||\}\leq
O(1)\mathcal{E}_{k}(f;Q)\leq O(1)\omega(r)
|f|_{\dot{C}^{k,\omega}_{q}(X)}.
\end{array}
$$
This immediately implies that
$$
\lim_{Q\to x}P_{Q}(f)(x)=\lim_{Q\to x}P_{Q}(x)=\widetilde f(x)
$$
and gives the required estimate of $|\widetilde f(x)-P_{Q}(f)(x)|$ by
the right-hand side of (\ref{e5.18}).
\end{proof}

Hereafter we assume for simplicity that
\begin{equation}\label{e5.19}
|f|_{\dot{C}^{k,\omega}_{q}(X)}=1.
\end{equation}
In particular, in this case
\begin{equation}\label{e5.20}
\mathcal{E}_{k}(f;Q)\leq\omega(r_{Q}),\ \ \ Q\in\mathcal{K}_{X}.
\end{equation}
\begin{lemma}\label{l5.6}
Let $Q\subset K$ be cubes from $\mathcal{K}_{X}$ of radii $r$ and $R$, respectively, $r<R\leq 2diam\ \! X$. Let $\widetilde K$ be the cube centered at $c_{K}$ of radius $2R$. Then it is true that
\begin{equation}\label{e5.21}
\mathcal{E}_{1}(f;Q)\leq O(1)\left\{r\!\int_{r}^{2R}
\frac{\omega(t)}{t^{2}}\ \!dt+\frac{r}{R}|||f;Q\cap \widetilde K|||\right\}.
\end{equation}
\end{lemma}
\begin{proof}
Choose $J\in\N$ from the condition
$$
R\leq 2^{J}r<2R
$$
and let $Q_{j}$ be the cubes centered at $c_{Q}$ and of radii $r_{j}:=2^{j}r$, $j=0,1,\dots J-1$, and $Q_{J}:=\widetilde K$, $r_{J}:=2R$.
Then $\{Q_{j}\}_{0\leq j\leq J}\subset\mathcal{K}_{X}$ is an increasing sequence of cubes.
We also set
$P_{j}:=P_{Q_{j}}$, $0\leq j\leq J$, see (\ref{e5.16}) for
the definition of $P_{Q}\in\mathcal{P}_{k-1}$. Under these notations we get
\begin{equation}\label{e5.22}
\mathcal{E}_{1}(f;Q)\leq\left\{\mathcal{E}_{1}(f-P_{Q};Q)+
\sum_{j=0}^{J-1}\mathcal{E}_{1}(P_{j+1}-P_{j};Q)+\mathcal{E}_{1}(P_{\widetilde K};Q)
\right\}.
\end{equation}
The first summand clearly equals
$$
\mathcal{E}_{k}(f;Q)\leq\omega(r)\leq O(1)\ \!
r \!\int_{r}^{2R}\frac{\omega(t)}{t^{2}}
\ \!dt
$$
as is required.

To estimate the remaining terms we use two inequalities whose proofs are
postponed to the end.
\begin{itemize}
\item[(A)]
Let $p$ be a polynomial of degree $k-1$ and $Q\in\mathcal{K}_{X}$ be a 
cube of radius $r$. Then
\begin{equation}\label{e5.23}
\mathcal{E}_{1}(p;Q)\leq O(1)r\max_{|\alpha|=1}|||D^{\alpha}p;Q\cap X|||.
\end{equation}
\item[(B)]
Let, in addition, $\widetilde Q\in\mathcal{K}_{X}$ be a cube of radius
$\widetilde r$ containing $Q$. Then 
\begin{equation}\label{e5.24}
\max_{|\alpha|=1}|||D^{\alpha}p;Q\cap X|||\leq O(1)\frac{1}{\widetilde r}
|||p:\widetilde Q\cap X|||.
\end{equation}
\end{itemize}
Using these inequalities to estimate the $j$-th term in (\ref{e5.22}) we get
$$
r^{-1}\mathcal{E}_{1}(P_{j+1}-P_{j};Q)\leq O(1)\frac{1}{r_{j}}
|||P_{j+1}-P_{j};Q_{j}\cap X|||.
$$
By the definitions of $s$-sets, $\omega$ and (\ref{e5.20}), 
the norm on the right-hand side 
is at most
$$
O(1)\frac{1}{r_{j}}\left(\mathcal{E}_{k}(f;Q_{j})+\left(
\frac{\mathcal{H}_{s}(Q_{j+1}\cap X)}
{\mathcal{H}_{s}(Q_{j}\cap X)}\right)^{1/q}\mathcal{E}_{k}(f;Q_{j+1})
\right)
\leq O(1)\frac{\omega(r_{j})}{r_{j}}.
$$
Moreover, by the definition of $r_{j}$ we get
$$
\frac{\omega(r_{j})}{r_{j}}\leq O(1)\int_{r_{j}}^{r_{j+1}}
\frac{\omega(t)}{t^{2}}\ \!dt,\ \ \ 0\leq j\leq J-1.
$$
Summing the finally obtained estimates over $j$ we then have
$$
\sum_{j=0}^{J-1}\mathcal{E}_{1}(P_{j+1}-P_{j};Q)\leq O(1)\ \!
r\ \!\!\int_{r}^{2R}
\frac{\omega(t)}{t^{2}}\ \!dt.
$$
Using now (\ref{e5.23}) and (\ref{e5.24}) we bound the last summand in
(\ref{e5.22}) by
$$
O(1)r\frac{|||P_{\widetilde K};\widetilde K\cap X|||}{R}\leq O(1)r\frac{2\ \! |||f;\widetilde K\cap X|||}{R}
$$
as is required.

To complete the proof of the lemma it remains to prove (\ref{e5.23}) and
(\ref{e5.24}). By the hypothesis of (A) we get
$$
\mathcal{E}_{1}(p;Q)\leq\inf_{\widetilde p}||p-\widetilde p||_{C(Q)}
\leq O(1)r\max_{|\alpha|=1}||D^{\alpha}p||_{C(Q)}
$$
where $\widetilde p$ runs over the space of polynomials of degree $0$.
The second of these inequalities is proved as follows. Using a 
homothety of $\Re^{n}$ we replace $Q$ by the unit cube 
$Q_{0}:=[0,1]^{n}$. The functions in $p$ of the both parts of this inequality
are norms on the finite-dimensional factor-space 
$\mathcal{P}_{k-1}/\mathcal{P}_{0}$ and therefore they are equivalent.
This implies the desired inequality.

Continuing the derivation we now use the weak Remez type inequality,
see (\ref{weak1}), and the fact that $X$ is an $s$-set to have
$$
||D^{\alpha}p||_{C(Q)}\leq O(1)|||D^{\alpha}p;Q\cap X|||
$$
and this completes the proof of (\ref{e5.23}).

Inequality (\ref{e5.24}) is proved in a similar way by means of the Markov inequality.
\end{proof}
\begin{lemma}\label{l5.6a}
$f=\widetilde f$ modulo $\mathcal{H}_{s}$-measure zero.
\end{lemma}
\begin{proof}
Let $L(f)$ be the Lebesgue set of $f$, i.e., the set of points $x\in X$ such that
$$
f(x)=\lim_{r\to 0}\frac{1}{\mathcal{H}_{s}(X_{r}(x))}\int_{X_{r}(x)}f\ \!d\mathcal{H}_{s}.
$$
Since $X$ is an $s$-set, the family of ``balls`` $\{X_{r}(x)\ :\ x\in X,\
0<r\leq 1\}$ satisfies axioms (i), (ii) in [St, p.8]. Therefore the
Corollary of Section I.3 from this book can be applied to our case with the measure $\mu:=\mathcal{H}_{s}|_{X}$. By this Corollary
$$
\mathcal{H}_{s}(X\setminus L(f))=0.
$$
It remains to show that
$$
f(x)=\widetilde f(x)\ \ \ {\rm for}\ \ \ x\in L(f).
$$
To this end  choose a cube $Q=Q_{r}(x)\in\mathcal{K}_{X}$, $0<r<1$, and set
$$
f_{r}(x):=\frac{1}{\mathcal{H}_{s}(X_{r}(x))}\int_{X_{r}(x)}f\ \!d\mathcal{H}_{s}.
$$
By the triangle inequality, the weak Remez type inequality for
$f_{r}(x)-P_{Q}$, see (\ref{weak1}), and the fact that $X$ is an $s$-set we obtain
\begin{equation}\label{prev}
|f_{r}(x)-P_{Q}(x)|\leq O(1)\{|||f-f_{r}(x);Q\cap X|||+\mathcal{E}_{k}(f,Q)\}.
\end{equation}
But $f\mapsto f_{r}$ is a projection from $L_{1}(X_{r}(x))$ onto the space $\mathcal{P}_{0}$ of polynomials of degree $0$ whose norm is $1$. Applying
an argument similar to that of Lemma \ref{l5.4} with this projection substituted for $\pi_{Q}$ we obtain that
$$
|||f-f_{r}(x);Q\cap X|||\leq O(1)\mathcal{E}_{1}(f;Q),
$$
and therefore by Lemma \ref{l5.6} and (\ref{e5.20}) for a sufficiently small $r$ the right-hand side of (\ref{prev}) is bounded by
$$
O(1)\{\mathcal{E}_{1}(f;Q)+\mathcal{E}_{k}(f;Q)\}\leq O(1)\left\{r\left(\int_{r}^{2}\frac{\omega(t)}{t^{2}}\ \!dt+|||f;K\cap X|||\right)+\omega(r)\right\}
$$
for some fixed cube $K$ of radius $1$ containing $Q$.  We conclude from here that for every $0<\epsilon<2$
$$
\begin{array}{c}
\displaystyle
\lim_{r\to 0}|f_{r}(x)-P_{Q}(x)|\leq\\
\\
\displaystyle
 O(1)\limsup_{r\to 0}\left(\omega(r)+r\left(\int_{r}^{\epsilon}\frac{\omega(t)}{t^{2}}\ \!dt+\int_{\epsilon}^{2}\frac{\omega(t)}{t^{2}}\ \!dt +|||f;K\cap X|||\right)\right)=\\
\\
\displaystyle
O(1)\limsup_{
r\to 0}\left(r\int_{r}^{\epsilon}\frac{\omega(t)}{t^{2}}\ \!dt\right)\leq O(1)\omega(\epsilon).
\end{array}
$$
Letting $\epsilon\to 0$ and noting that $\lim_{r\to 0}f_{r}(x)=f(x)$ for the Lebesgue point $x$ and $\lim_{Q\to x}P_{Q}(x)=\widetilde f(x)$ we complete the proof of the lemma.
\end{proof}

Now we finalize the proof of Theorem \ref{ybr}. For $Q\in\mathcal{K}_{X}$ and the polynomial $P_{Q}(f)$ of degree $k-1$ defined in (\ref{e5.12}) we set
$$
\widetilde P_{Q}(f):=P_{Q}(f)-P_{Q}(f)(c_{Q})+\widetilde f(c_{Q}).
$$
Then $\widetilde P_{Q}(f)(c_{Q})=\widetilde f(c_{Q})$ and Claim 2, see (\ref{e5.14}), is true for the family $\widetilde\Pi(f):=\{\widetilde P_{Q}\}_{Q\in\mathcal{K}_{X}}$. Show that Claim 1 is also true for $\widetilde\Pi(f)$.

Let $Q\subset Q'$ be cubes from $\mathcal{K}_{X}$ of radii $r<r'$ satisfying for some $i$ the condition
$$
t_{i}\leq r<r'\leq t_{i+2}.
$$
By the weak Remez type inequality, see (\ref{weak1}), and Lemma \ref{l5.4} we have
$$
\begin{array}{c}
\displaystyle
\max_{Q}|P_{Q}(f)-P_{Q'}(f)|\leq O(1)\max_{X\cap Q}|P_{Q}(f)-P_{Q'}(f)|\leq
\\
\\
\displaystyle
O(1)|||P_{Q}(f)-P_{Q'}(f);X\cap Q|||\leq
O(1)\!\left\{\mathcal{E}_{k}(f;Q)+\left(\frac{\mathcal{H}_{s}(Q'\cap X)}{\mathcal{H}_{s}(Q\cap X)}\right)^{1/q}\mathcal{E}_{k}(f;Q')\right\}.
\end{array}
$$
Both of the best approximations are bounded by $\omega(r')|f|_{\dot{C}_{q}^{k,\omega}(X)}$ while, since $X$ is an $s$-set, the ratio of $\mathcal{H}_{s}$-measures is at most
$$
\left\{\frac{a}{b}\left(\frac{r'}{r}\right)^{s}\right\}^{1/q}\leq O(1)\left(\frac{t_{i+2}}{t_{i}}\right)^{s/q}\leq O(1).
$$
Hence, in this situation, see (\ref{e5.19}),
$$
\max_{Q}|P_{Q}(f)-P_{Q'}(f)|\leq O(1)\omega(r').
$$
Moreover, by Lemma \ref{l5.5a}, see(\ref{e5.18}),
$$
|\widetilde f(c_{Q})-P_{Q}(f)(c_{Q})|\leq O(1)\omega(r).
$$
Taking into account the definition of $\widetilde P_{Q}(f)$ we then obtain the inequality
$$
\max_{Q}|\widetilde P_{Q}(f)-\widetilde P_{Q'}(f)|\leq O(1)\omega(r')
$$
as is required in the definition of a $(k,\omega,X)$-chain.

This completes the proof of Claim 1 and therefore of Theorem \ref{ybr}. 
\end{proof}


\begin{thebibliography}{ }
\bibitem[ABr]{Br}
A. Brudnyi, On a BMO-property for subharmonic functions. J. Fourier
Anal. Appl. {\bf 8} (6) (2002), 603-612.
\bibitem[BB]{BB}
A. Brudnyi and Yu. Brudnyi, Local inequalities for multivariate
polynomials and plurisubharmonic functions. In Frontiers in
Interpolation and Approximation. Eds. Govil et al. Chapman $\&$
Hall, 2006, 17-32.
\bibitem[BG]{BG}
Yu. A. Brudnyi and M. I. Ganzburg, On an extremal problem for polynomials
of $n$-variables. Math. USSR Izv. {\bf 37} (1973), 344-355.
\bibitem[BSh1]{BSh1}
Yu. A. Brudnyi and P. A. Shvartsman, Description of trace of general
Lipschitz space to an arbitrary compact. (Russian) 
In: Studies in the theory of
functions of several real variables. Yaroslavl State Univ. Yaroslavl, 1982,
16-24.
\bibitem[BSh2]{BSh2}
Yu. A. Brudnyi and P. A. Shvartsman, The Whitney problem of existence of
a linear extension operator. J. Geom. Anal. {\bf 7} (1997), no. 4, 515-574.
\bibitem[Ca]{Ca}
S. Campanato, Proprieta di una famiglia di spazi funzional. Ann.
Scuola Norm. Sup. Pisa {\bf 18} (1964), 137-160.
\bibitem[C]{C}
H. Cartan, Sur les syst\`{e}mes de fonctions holomorphes \'{a}
vari\'{e}t\'{e}s lin\'{e}aires lacunaires et leurs applications.
Ann. Sci. \'{E}cole Norm. Sup. {\bf 45} (3) (1928), 255-346.
\bibitem[F]{F}
H. Federer, Geometric measure theory. Springer-Verlag. New York,
1969.
\bibitem[G]{G}
M. I. Ganzburg, Polynomial inequalities on measurable sets and their
applications. Constr. Approx. {\bf 17} (2001), 275-306.
\bibitem[GH]{GH}
P. Griffiths and J. Harris, Principles of Algebraic Geometry. New
York, 1978.
\bibitem[GK]{GK}
E. A. Gorin and A. L. Koldobskii, On potentials of measures with
values in a Banach space. Sibirsk. Mat. Zh. {\bf 28} (1), 65-80.
English transl. in Siberian Math. J. {\bf 28} (1987).
\bibitem[J]{J}
P. W. Jones, Quasiconformal mappings and extendability of functions
in Sobolev classes. Acta Math. {\bf 147} (1981), 71-88.
\bibitem[JN]{JN}
F. John and L. Nirenberg, On functions of bounded mean oscillation.
Comm. Pure Appl. Math. {\bf 14} (1961), 415-426.
\bibitem[JW]{JW}
A. Jonsson and H. Wallin, Functions spaces on subsets of
$\mathbf{R}^{n}$. Harwood Academic Publishers, 1984.
\bibitem[KS]{KS}
M. I. Kadets and M. G. Snobar, Certain functionals on the Minkowski 
compactum. (Russian) Mat. Zametki {\bf 10} (1971), 453--457; English
transl. in Math. Notes {\bf 10} (1971).
\bibitem[Ma]{Ma}
P. Mattila, Geometry of sets and measures in Euclidean spaces.
Cambridge Univ. Press, 1995.
\bibitem[M]{M}
D. Mumford, Algebraic Geometry I. Complex projective varieties.
Springer, 1976.
\bibitem[Mo]{Mo}
C. B. Morrey, Multiple integrals in the calculus of variations. Springer, 
1966.
\bibitem[L]{L}
B. Ya. Levin, Lectures on entire functions. Amer. Math. Soc. Transl.
of Math. Monographs {\bf 150}, Providence RI, 1996.
\bibitem[St]{St}
E. M. Stein, Singular integrals and differential properties of functions.
Princeton Univ. Press. Princeton, 1970.
\bibitem[Tr]{Tr}
H. Triebel, Fractals and spectra related to Fourier analysis and
function spaces. Monographs in Math. vol 91. Birkh\"{a}user Verlag,
1997.
\bibitem[YBr]{YBr}
Yu. A. Brudnyi, On an extension theorem. Funk. Anal. i Prilozhen. {\bf
4} (1970), 97-98; English transl. in Funct. Anal. Appl. {\bf 4}
(1970), 252-253.

\end{thebibliography}
\end{document}